\documentclass[11pt]{article}
\usepackage{amsfonts}
\usepackage{latexsym}
\usepackage[all]{xy}

\usepackage{amsfonts}

\begin{document}

\author{Hans-Peter A. K\"unzi, S. Romaguera and M.A. S\'anchez Granero}
\title{The bicompletion of the Hausdorff
quasi-uniformity}
\date{{\small Dedicated to Professor Robert Lowen on the occasion of his sixtieth birthday}
\\ \bigskip 30 September, 2008
}

\newtheorem{Prp}{\qquad Proposition}
\newtheorem{Thm}{\qquad Theorem}
\newtheorem{Lem}{\qquad Lemma}
\newtheorem{Cor}{\qquad Corollary}
\newtheorem{Rem}{\qquad Remark}
\newtheorem{Ex}{\qquad Example}
\newtheorem{Prb}{\qquad Problem}
\newtheorem{Def}{\qquad Definition}
\maketitle

\begin{abstract}
We study conditions under which the Hausdorff
quasi-uniformity ${\mathcal U}_H$ of a quasi-uniform
space $(X,{\mathcal U})$ on the set ${\mathcal P}_0(X)$
of the nonempty
subsets of $X$ is bicomplete.

Indeed we present an explicit method to construct
the bicompletion of the $T_0$-quotient of the Hausdorff quasi-uniformity of a quasi-uniform space.
It is used to find a characterization
of those quasi-uniform $T_0$-spaces $(X,{\mathcal U})$ for which the Hausdorff quasi-uniformity
$\widetilde{{\mathcal U}}_H$ of their bicompletion $(\widetilde{X},{\widetilde{\mathcal U}})$
on ${\mathcal P}_0(\widetilde{X})$ is bicomplete.

\noindent
\footnotetext{\noindent
AMS (2000) Subject Classifications: 54E15 54B20 54D35 54E55

}

\noindent \footnotetext{The first author would like
to thank the South African Research Foundation for partial
financial support under grant FA2006022300009. He also thanks the Polytechnical
University of Valencia for its hospitality during his stay at this institution during
the summer of 2007.
The second and third authors thank the support of the Spanish Ministry of
Education and Science and FEDER, grant MTM2006-14925-C02-01.}

\end{abstract}

\section{Introduction}
In the theory of quasi-uniform spaces the construction
of the bicompletion is well known
(see e.g. \cite[Theorem 3.33]{FletcherLindgren}).
Also the Hausdorff
quasi-uniformity of a quasi-uniform space was investigated by many authors
(see e.g.
\cite{Berthiaume,CaoReillyRomaguera,KunziRomaguera,KunziRyser,
RodriguezRomaguera,SanchezGranero}).
In this article we want to study the problem under which
conditions the Hausdorff quasi-uniformity ${\mathcal U}_H$
of a quasi-uniform space $(X,{\mathcal U})$
on the set ${\mathcal P}_0(X)$ of nonempty subsets of $X$
is bicomplete. Some
 results dealing with our question can be found
 in the article of K\"unzi and Ryser \cite{KunziRyser}.
In particular these authors observed that
 the Hausdorff quasi-uniformity ${\mathcal U}_H$
 of a
   totally bounded and bicomplete quasi-uniformity
   ${\mathcal U}$ is (totally bounded and)
bicomplete  (see \cite[Corollary 9]{KunziRyser}). Recall that a
quasi-uniform space $(X,{\mathcal U})$ is totally bounded and
bicomplete if and only if the topology $\tau({\mathcal U}^s)$ is
compact, where ${\mathcal U}^s$ denotes the coarsest uniformity
finer than ${\mathcal U}$ (see e.g. \cite[Proposition
2.6.10]{KunziIntroduction}).

We also note that our question has a well-known and
satisfactory answer in the setting of uniform
spaces. To this end recall the Burdick-Isbell
\cite{Burdick,Isbell} result which says that for a {\em uniform}
space $(X,{\mathcal U})$ the Hausdorff uniformity ${\mathcal U}_H$ on ${\mathcal P}_0(X)$ is
complete if and only if each stable
filter on $(X,{\mathcal U})$ has a cluster point. For uniform spaces the latter property is usually
called {\em supercompleteness} and has been investigated by
many authors (see e.g. \cite{ArticoMarconiPelant,BuhagiarPasynkov,FedorchukKunzi,Hohti}).
It is well known that a uniform space with a countable
base
is complete if and only if it is supercomplete
 (see for instance the discussion preceding \cite[Lemma 3.4.7]{KunziIntroduction}).
For an application of the Burdick-Isbell condition to the theory
of topological groups we refer the reader to
\cite{RomagueraSanchis}. A quasi-uniform variant of the
Burdick-Isbell result was obtained by K\"unzi and Ryser
\cite[Proposition 6]{KunziRyser} who proved that for a
quasi-uniform space $(X,{\mathcal U})$ the Hausdorff
quasi-uniformity ${\mathcal U}_H$ on ${\mathcal P}_0(X)$ is right
$K$-complete if and only if each stable filter on $(X,{\mathcal
U})$ has a cluster point. Further related investigations about
quasi-uniform spaces were conducted by S\'anchez-Granero
\cite{SanchezGranero} and Burdick \cite{Burdick2}.

In this article we first discuss a characterization of
bicompleteness of the Hausdorff quasi-uniformity due to K\"unzi
and Ryser. Then we present a general method to construct the
bicompletion of the $T_0$-quotient of the Hausdorff
quasi-uniformity of a quasi-uniform space. The result is finally
used to find a condition under which for a quasi-uniform
$T_0$-space $(X,{\mathcal U})$ the Hausdorff quasi-uniformity
$\widetilde{{\mathcal U}}_H$ on ${\mathcal P}_0(\widetilde{X})$ of
the bicompletion $(\widetilde{X},\widetilde{\mathcal U})$ of
$({X},{\mathcal U})$ is bicomplete.

\bigskip
For the basic facts about quasi-uniformities we refer the
reader to \cite{FletcherLindgren} and \cite{KunziIntroduction}.
In particular for a quasi-uniform space $(X,{\mathcal U})$ the filter
${\mathcal U}^{-1}$ on $X\times X$  denotes the conjugate quasi-uniformity
of ${\mathcal U},$ and ${\mathcal U}^s$, as already mentioned above, denotes the coarsest uniformity finer
than ${\mathcal U}$ on $X.$ Similarly, for an entourage $U$ of a quasi-uniform
space $(X,{\mathcal U}),$ $U^s$ denotes the relation $U\cap U^{-1}.$

We recall that a quasi-uniform space $(X,{\mathcal U})$
is called {\em bicomplete} provided that the uniformity
${\mathcal U}^s$ is complete.
It is well known that a quasi-pseudometric $d$ is
(sequentially) bicomplete if and only if the induced
quasi-pseudometric quasi-uniformity
${\mathcal U}_d$ is bicomplete (see e.g. \cite[beginning of Section 2.6]{KunziIntroduction}).
Each quasi-uniform $T_0$-space $(X,{\mathcal U})$
can be embedded into a(n up-to quasi-uniform isomorphism unique) bicomplete quasi-uniform $T_0$-space
$(\widetilde{X},\widetilde{\mathcal U})$ (its so-called {\em bicompletion}) in which it is
$\tau({\widetilde{\mathcal U}^s})$-dense \cite[Theorem 3.33]{FletcherLindgren}.
An explicit construction of the bicompletion $(\widetilde{X},\widetilde{\mathcal U})$
of a quasi-uniform $T_0$-space $(X,{\mathcal U})$ is described below.

Given a quasi-uniform space $(X,{\mathcal U})$ we shall consider the
${\mathcal U}$-{\em equivalence relation} $\sim$ on $X$ that underlies its
$T_0$-reflection: For
$x,y\in X$ we
have  $x\sim y$ if and only if $(x,y)\in \bigcap {\mathcal U}\cap (\bigcap
{\mathcal U}^{-1}).$
It is well known that the $T_0$-quotient of a quasi-uniform space $(X,{\mathcal U})$
can be represented by any subspace of $(X,{\mathcal U})$ that intersects
each ${\mathcal U}$-equivalence
class exactly in a singleton.

By $\mbox{adh}_\tau{\mathcal F}$
we shall denote the set of cluster points of a filter ${\mathcal F}$ on a topological
space $(X,\tau)$
with respect
to the topology $\tau.$
A filter ${\mathcal F}$ on a quasi-uniform space
$(X,{\mathcal U})$ is called {\em stable}
provided that $\bigcap_{F\in {\mathcal F}}U(F)\in
{\mathcal F}$ whenever $U\in {\mathcal U}.$

\begin{Lem} \label{subs} Let $A$ be a subset of a quasi-uniform
space $(X,{\mathcal U}).$ If ${\mathcal F}$
is a stable filter on the subspace
$(A,{\mathcal U}\vert A)$
of $(X,{\mathcal U}),$ then the filter generated
by the filterbase ${\mathcal F}$ on $(X,{\mathcal U})$
is stable on $(X,{\mathcal U}).$
\end{Lem}

{\em Proof.}  The assertion is obvious. $\hfill{\Box}$
\bigskip

Recall finally that a quasi-uniform space $(X,{\mathcal U})$
is called {\em precompact} provided that for each
$U\in {\mathcal U}$ there is a finite subset $F$
of $X$ such that $\bigcup_{x\in F}U(x)=X.$
A quasi-uniform space $(X,{\mathcal U})$ is said to be  {\em totally bounded}
provided that ${\mathcal U}^s$ is precompact.
Totally bounded quasi-uniform spaces are precompact, but the converse does not hold.

\section{Preliminaries}

In order to discuss the investigations of K\"unzi
and Ryser \cite{KunziRyser} that are relavant to our problem it is useful to recall
first several additional concepts.
Let $(X,\rho,\sigma)$ be a bitopological space.
The {\em double closure} of a set $C\subseteq X$ is
$\mbox{cl}_{\rho}C\cap \mbox{cl}_{\sigma}C.$
A subset of a bitopological space $(X,\rho,\sigma)$ is called {\em
doubly closed}
if it is equal to its double closure.
Observe
that each $\rho$-closed set as well as each $\sigma$-closed set
is doubly closed.

Of course, the intersection of an arbitrary family of doubly closed
sets is doubly closed.
Indeed the double closure operator is an (idempotent) closure
 operator, which
in general does not commute with finite unions. Therefore it is
not a (topological) Kuratowski closure operator. For instance the
union of two doubly closed sets need not be doubly closed: The
intervals $]0,1[$ and $]1,2[$ are both doubly closed\footnote{If
we speak about doubly closed subsets of a quasi-uniform space
$(X,{\mathcal U)},$ then we always mean doubly closed with respect
to the bitopological space $(X,\tau({\mathcal U}),\tau({\mathcal
U}^{-1})).$} in the Sorgenfrey line $(\mathbb{R},{\mathcal U}_s)$
(for a definition of this quasi-uniform space see Section
\ref{Sorgen}). But $1$ clearly belongs to the double closure of
the union of these intervals.

In a quasi-uniform $T_0$-space $(X,{\mathcal U})$
each singleton is doubly closed.
For a quasi-uniform space $(X,{\mathcal U}),$
the nonempty doubly closed subsets of $X$ can represent the
${\mathcal U}_H$-equivalence
classes of the space
$({\mathcal P}_0(X),{\mathcal U}_H)$ (see \cite[Lemma 2]{KunziRyser}).
Indeed each nonempty subset $A$ of $X$
is ${\mathcal U}_H$-equivalent to its double closure
$\mbox{cl}_{\tau({\mathcal U})}A\cap
\mbox{cl}_{\tau({\mathcal U}^{-1})}A.$

By definition the set $C$ of
{\em double cluster points} of a filter
${\mathcal F}$ on
a bitopological space $(X,\rho,\sigma)$
 is the adherence
of ${\mathcal F}$ with respect to the
double closure operator. Hence $C=
\bigcap_{F\in {\mathcal F}} (\mbox{cl}_\rho
F\cap \mbox{cl}_\sigma F)= \mbox{adh}_\rho{\mathcal F}\cap \mbox{adh}_\sigma{\mathcal F}.$

Let ${\mathcal F}$ be a filterbase on a quasi-uniform space
$(X,{\mathcal U}).$
For each $U\in {\mathcal U}$ set
 $U_{\mathcal F}:=\bigcap_{F\in {\mathcal
F}}(U^{-1}(F)\cap
   U(F)).$
A filterbase ${\mathcal F}$ on a quasi-uniform
space $(X,{\mathcal U})$ is called {\em doubly stable}
provided that $U_{\mathcal F}$
   belongs to ${\mathcal F}$
   whenever $U\in {\mathcal U}.$
Hence a filter on a quasi-uniform space $(X,{\mathcal U})$ is doubly
stable if and only if it is ${\mathcal U}$-stable and ${\mathcal U}^{-1}$-stable.
For a doubly stable filter ${\mathcal F}$
on a quasi-uniform space $(X,{\mathcal U})$
we shall denote the filter on $X$ generated by the filterbase
$\{U_{\mathcal F}:U\in {\mathcal U}\}$
by ${\mathcal F}_{\mathcal U}.$ Of course, we have
${\mathcal F}_{\mathcal U}\subseteq {\mathcal F}.$

Note that each ${\mathcal U}^s$-Cauchy filter on a quasi-uniform
space $(X,{\mathcal U})$ is ${\mathcal U}^s$-stable, and thus
doubly stable. These three concepts for filters on a quasi-uniform
space coincide for ultrafilters, but not in general (compare
\cite[Proposition 2.6.5]{KunziIntroduction}). For instance the
cofinite filter on the (bicomplete) Sorgenfrey line
$(\mathbb{R},{\mathcal U}_s)$
 is a doubly stable filter, which is not $({\mathcal
 U}_s)^s$-stable and
does not have a $\tau(({\mathcal U}_s)^s)$-limit point. Furthermore each real number is a double cluster point
of this filter. Note also that each complete uniform space which is not supercomplete
has a stable filter without cluster point.

\begin{Lem} A quasi-uniform space $(X,{\mathcal U})$ is totally bounded if and only if
each filter on $(X,{\mathcal U})$
is doubly stable.
\end{Lem}

{\em Proof.} This observation follows from the following two
results: A quasi-uniform space $(X,{\mathcal U})$ is totally
bounded if and only if both ${\mathcal U}$ and ${\mathcal U}^{-1}$
are hereditarily precompact \cite[Corollary
9]{KunziMrsevicReillyVamanamurthy}. For a quasi-uniform space
$(X,{\mathcal U})$ the quasi-uniformity ${\mathcal U}^{-1}$ is
hereditarily precompact if and only if each filter on $X$ is
${\mathcal U}$-stable \cite[Proposition 2.5]{KunziJunnila}.
$\hfill{\Box}$

\bigskip
   The proof of the following remark is also immediate.

\begin{Rem} Given
a quasi-uniform space $(X,{\mathcal U}),$
for any $C\in {\mathcal P}_0(X)$
the filter ${\mathcal C}_C$ on $X$
generated by the filterbase $\{C\}$ is doubly stable; in fact it is
${\mathcal U}^s$-stable.
$\hfill{\Box}$
\end{Rem}

In  \cite[Proposition 8]{KunziRyser} K\"unzi and Ryser
 showed that the Hausdorff quasi-uniformity
${\mathcal U}_H$
of a quasi-uniform space $(X,{\mathcal U})$
is bicomplete if and only if each doubly stable filter ${\mathcal F}$
on $(X,{\mathcal U})$ satisfies the following condition:

For any $U\in {\mathcal U}$ there is $F\in {\mathcal F}$ such that
$F\subseteq U^{-1}(C({\mathcal F}))\cap U(C({\mathcal F})),$ where
$C({\mathcal F})$ denotes the set of double cluster points of
${\mathcal F}.$ (Note that in particular this condition implies
that each doubly stable filter has a double cluster point. Since a
${\mathcal U}^s$-Cauchy filter on a quasi-uniform space
$(X,{\mathcal U})$ that has a double cluster point $\tau({\mathcal
U}^s)$-converges to this point (compare e.g. \cite[Proposition
2.6.5]{KunziIntroduction}), it immediately follows  that a
quasi-uniform space $(X,{\mathcal U})$ is bicomplete if
$({\mathcal P}_0(X),{\mathcal U}_H)$ is bicomplete.)

\medskip

In the following we shall call the aforementioned
condition the {\em K\"unzi-Ryser condition}.
That condition for
a quasi-uniform space $(X,{\mathcal U})$ can be reformulated in a way that reveals
how it is
related to $\tau({\mathcal U}^s)$-compactness of a quasi-uniform space
$(X,{\mathcal U}),$ that is,
the property that each filter on $X$
has a $\tau({\mathcal U}^s)$-cluster point:
For each doubly stable filter ${\mathcal F}$ on $(X,{\mathcal U})$ and $U\in {\mathcal U}$
there exists $V\in {\mathcal U}$ so that
for all $x\in X$ such that ${\mathcal F}$ traces on $V(x)$ and
$V^{-1}(x),$\footnote{A filter ${\mathcal F}$ on a set $X$
{\em traces on a subset $A$ of $X$} provided that
$F\cap A\not =\emptyset$ whenever $F\in {\mathcal F}.$} both $U(x)$ and $U^{-1}(x)$ contain a double cluster point of
${\mathcal F}$ (which may be distinct).
Indeed that property can be stated in the following way, which
shows that it is equivalent to the K\"unzi-Ryser condition:
For each doubly stable filter ${\mathcal F}$ on $(X,{\mathcal U})$ and each $U\in {\mathcal U}$
there is $V\in {\mathcal U}$ such that
$V_{\mathcal F}\subseteq U^{-1}(C({\mathcal F}))\cap
   U(C({\mathcal F}))$ where $C({\mathcal F})$ is the set of double cluster points of
   ${\mathcal F}.$

We next  present a simple example illustrating
the K\"unzi-Ryser condition.

\begin{Ex} \label{bei} Let $(X,{\mathcal U})$ be a quasi-uniform space possessing
  some entourage $V\in {\mathcal U}$ such that for each $x\in X,$
$V(x)=\{x\}$ or $V^{-1}(x)=\{x\}.$
Then $({\mathcal P}_0(X),{\mathcal U}_H)$ is bicomplete.
\end{Ex}

{\em Proof.} Let ${\mathcal F}$ be a doubly stable filter
on $(X,{\mathcal U}).$ Then $\bigcap_{F\in {\mathcal F}}(V^{-1}(F)
\cap V(F))\in {\mathcal F}.$
Let $x\in \bigcap_{F\in {\mathcal F}}(V^{-1}(F)
\cap V(F)).$ By assumption $V^{-1}(x)=\{x\}$
or $V(x)=\{x\}.$ Thus in either case $x\in
\bigcap_{F\in {\mathcal F}}F\subseteq
\mbox{adh}_{\tau({\mathcal U})}{\mathcal F}\cap
\mbox{adh}_{\tau({\mathcal U}^{-1})}{\mathcal F}
\subseteq C({\mathcal F}).$
It follows that $\bigcap_{F\in {\mathcal F}}(V^{-1}(F)\cap
V(F))=C({\mathcal F}).$
We conclude  that the K\"unzi-Ryser condition is satisfied,
since $C({\mathcal F})\in {\mathcal F}.$
$\hfill{\Box}$

\begin{Rem} Note that the usual quasi-uniformity
\cite[Corollary 7]{KunziRyser} on the Mr\'owka space
$\Psi$ satisfies
the condition of Example \ref{bei}.
\end{Rem}

\section{A counterexample to a possible weakening of the
K\"unzi-Ryser condition}

In the light of the aforementioned Burdick-Isbell condition
that characterizes completeness of the Hausdorff uniformity
of a uniform space it is natural to conjecture that the K\"unzi-Ryser condition
is unnecessarily complicated and that
the Hausdorff quasi-uniformity of a quasi-uniform
space $(X,{\mathcal U})$ is bicomplete if and only if each doubly stable filter on
$(X,{\mathcal U})$ has a double
cluster point. However the following quasi-uniform space yields a counterexample
to that conjecture.

\begin{Ex} \label{contra} There exists a quasi-pseudometrizable quasi-uniform space
$(X,{\mathcal U})$
such that each doubly stable filter on $(X,{\mathcal U})$ has a double
cluster point,
although the K\"unzi-Ryser condition is not satisfied. (Hence the
Hausdorff quasi-uniformity
${\mathcal U}_H$ of ${\mathcal U}$ is not bicomplete.)
\end{Ex}

{\em Proof.} Let $X$ be the set $\mathbb{N}$ of positive integers and let ${\mathcal U}$ be the
quasi-uniformity
generated by the countable subbase consisting of the usual order $\leq$ on $X$ and all the transitive relations
$[\{x\}\times X]\cup [X\times (X\setminus \{x\})]$ whenever $x\in X.$
Note that $\tau({\mathcal U})$ is the cofinite topology on $X$ and
 $\tau({\mathcal U}^{-1})$
is the discrete topology on $X.$

We first show that each stable filter on $(X,{\mathcal U})$ has a
$\tau({\mathcal U}^s)$-cluster point.
Therefore in particular each doubly stable filter on $(X,{\mathcal U})$
has a double cluster point
in $(X,{\mathcal U}).$
Indeed let ${\mathcal F}$ be a stable filter on $(X,{\mathcal U}).$
If $\bigcap {\mathcal F}$ is
nonempty, then any point in that intersection is clearly a
 $\tau({\mathcal U}^s)$-cluster
point and
we are finished.
So we can assume that $\bigcap {\mathcal F}=\emptyset.$
Then for each $n\in X$ there is $F_n\in {\mathcal F}$ such that $F_n$ does
not contain any positive integer smaller than $n.$ We conclude that for $U=\leq$
we have that $\bigcap_{F\in {\mathcal F}}U(F)=\emptyset.$ Hence we have
reached a contradiction, since it follows that ${\mathcal F}$
is not stable on $(X,{\mathcal U}).$ Therefore for each stable filter
${\mathcal F}$
on $(X,{\mathcal U})$
we have $\bigcap_{F\in {\mathcal F}}F \not =
\emptyset$ and we are finished.

Next we consider the filter ${\mathcal G}$ on $X$ generated by the
base $\{\{1\}\cup G':$
   $G'$ is a cofinite subset of $X\}.$
First note that for any $U\in {\mathcal U}$ and $G\in {\mathcal G}$
   we have $U^{-1}(G)=X,$ because $G$ is cofinite and for each $y\in X,$ $U(y)$
is cofinite and thus $G$ and $U(y)$ intersect.
In particular we conclude that
${\mathcal G}$ is ${\mathcal U}^{-1}$-stable.

Furthermore, for any $U\in {\mathcal U}$ we have $U(1)\subseteq U(G)$
whenever $G\in
{\mathcal G}.$ Since $U(1)$ is a cofinite set that contains $1,$ it belongs
to the filter ${\mathcal G}.$ Hence we have shown that ${\mathcal G}$
is ${\mathcal U}$-stable.
Consequently ${\mathcal G}$ is a doubly stable filter on $(X,{\mathcal U}).$

Clearly $\{1\}=\bigcap_{G\in {\mathcal G}}G.$
We conclude that $1$ is a double cluster point
of ${\mathcal G},$ and the only one of ${\mathcal G}$, because the topology
$\tau({\mathcal U}^{-1})$ is discrete.
Therefore $C({\mathcal G}):=\{1\}$ is the set of double cluster points of ${\mathcal G}.$
Set $U=\leq.$
We see that $U^{-1}(1)=\{1\}.$ It follows that
   $G\not \subseteq U^{-1}(C({\mathcal G}))$ whenever $G\in {\mathcal G}.$
Thus ${\mathcal G}$ is a doubly stable filter on $(X,{\mathcal U})$
that does not satisfy the K\"unzi-Ryser condition.
Consequently ${\mathcal U}_H$ is not  bicomplete.

It is also interesting to note that the filter generated by
$\{G\setminus U^{-1}(C({\mathcal G})):G\in {\mathcal G}\}$ is equal to the cofinite
filter on $X$ and therefore not contained in a
stable filter on $(X, {\mathcal U})$ (compare with \cite[Proof of Lemma 6]{KunziRyser}
where a similar construction for filters on a quasi-uniform space is studied
 that preserves stability  of filters).

Our bicomplete example $(X,{\mathcal U})$
also has the property that each (doubly) stable filter  is contained in a
$\tau({\mathcal U}^s)$-neighborhood filter,
but nevertheless the Hausdorff quasi-uniformity ${\mathcal U}_H$ is
not bicomplete on ${\mathcal P}_0(X)$ (compare Proposition \ref{unifor} below).
$\hfill{\Box}$

\begin{Rem} We recall that a quasi-uniform space $(X,{\mathcal U})$ is called
{\em half-complete} provided that each ${\mathcal U}^s$-Cauchy filter
$\tau({\mathcal U})$-converges. In \cite[Proposition 3.13]{SanchezGranero}
those quasi-uniform spaces $(X,{\mathcal U})$ were characterized for which
$({\mathcal P}_0(X),{\mathcal U}_H)$ is half-complete.
With the help of this criterion one readily
 checks that the quasi-uniformity ${\mathcal U}_{H}$ in Example \ref{contra}
is half-complete, because $\tau({\mathcal U})$ is compact.
Together with the argument presented above about the filter
${\mathcal G}$ the criterion also establishes that $({\mathcal
U}^{-1})_H=({\mathcal U}_H)^{-1}$ is not half-complete, although
each doubly stable filter on $(X,{\mathcal U}^{-1})$ has a
$\tau({\mathcal U}^{-1})$-cluster point.$\hfill{\Box}$

\end{Rem}

\section{Another positive application of the K\"unzi-Ryser condition} \label{Sorgen}

Let $\mathbb{R}$ denote the set of the reals. As usual (see e.g.
\cite[Corollary 6]{KunziRyser}) define the Sorgenfrey quasi-metric
$s$ on $\mathbb{R}$ as follows: For each $x,y\in \mathbb{R}$  set
$s(x,y)=y-x$ if $y\geq x$ and $s(x,y)=1$ otherwise. In
\cite[Example 7]{KunziRyser} it was shown that a doubly stable
filter on the set $\mathbb{Q}$ of the rationals (equipped with the
(bicomplete) restriction of the quasi-uniformity ${\mathcal U}_s$
induced by
 the Sorgenfrey quasi-metric $s$) need not have a double cluster point in $\mathbb{Q}$.
Hence the corresponding Hausdorff
quasi-uniformity on ${\mathcal P}_{0}(\mathbb{Q})$ is not bicomplete.
In this section we are going to show that $\mathbb{R}$ behaves differently.

\begin{Ex}
The quasi-pseudometrizable quasi-uniform space
 $({\mathcal P}_0(\mathbb{R}),{({\mathcal U}_s)}_H)$ is bicomplete,
where $s$ denotes the Sorgenfrey quasi-metric on the set $\mathbb{R}$ of the reals.
\end{Ex}

{\em Proof.} For each $n\in \mathbb{N}$ set $S_{2^{-n}}=\{(x,y)\in \mathbb{R}\times \mathbb{R}: s(x,y)<2^{-n}\}.$
Let ${\mathcal F}$ be a doubly stable filter on $(\mathbb{R},{\mathcal U}_s).$
For each $n\in \mathbb{N},$ set $F_n=
\bigcap_{F\in {\mathcal F}}S_{2^{-n}}^{-1}(F)\cap \bigcap_
{F\in {\mathcal F}}S_{2^{-n}}(F).$ Observe that the sequence $(F_n)_{n\in \mathbb{N}}$ is
decreasing and
that the filter on $\mathbb{R}$ generated by the filterbase $\{F_n:n\in \mathbb{N}\}$ has the
same sets of cluster points with respect to the topologies $\tau(s)$ and $\tau(s^{-1})$ as
${\mathcal F}$ has.
Furthermore by assumption on ${\mathcal F},$ $F_n\in {\mathcal F}$ whenever $n\in \mathbb{N}$.
Note also that $\bigcap_{n\in \mathbb{N}} F_n
=C$ is the set of double cluster points of ${\mathcal F}$ in $(\mathbb{R},s),$ that is,
the set $\mbox{adh}_{\tau(s)}{\mathcal F}\cap
\mbox{adh}_{\tau(s^{-1})}{\mathcal F}.$

Consider an arbitrary $a\in \mathbb{N}.$ We show that the
assumption that $F_n\setminus S_{2^{-a+3}}(C)\not =\emptyset$
whenever $n\in \mathbb{N}$ leads to a contradiction. For each
$n\in \mathbb{N}$ set $E_n:=F_n\setminus S_{2^{-a+1}}(C).$ Let
${\mathcal E}$ be the filter generated by the filterbase
$\{E_n:n\in \mathbb{N}\}$ on $\mathbb{R}.$ Choose $x_a\in
F_a\setminus S_{2^{-a+3}}(C).$ Observe that
   $]x_a-2^{-a+2},x_a]=S^{-1}_{2^{-a+2}}(x_a)$ is disjoint from
$S_{2^{-a+1}}(C).$

   We first show that $\mbox{adh}_{\tau(s)}{\mathcal E}\cap \,
]x_a-2^{-a+2},x_a[$ is nonempty.
Given $x_n$ with $n\in \mathbb{N}$ and $n\geq a,$ by definition of $F_n$ inductively we
find $x_{n+1}\in
   F_{n+1}$ such that $x_n\in S_{2^{-n}}(x_{n+1}).$
Thus $s(x_{n+1},x_n)<2^{-n}$ whenever $n\in \mathbb{N}$ and $n\geq a.$
Note that $x_n\in E_n$ whenever $n\in \mathbb{N}$ and $n \geq a,$ since $x_n\in
S^{-1}_{2^{-a+2}}(x_a).$
Also the sequence $(x_n)_{n\geq a}$  converges to its infimum $x$
with respect to the topology $\tau(s).$
These results follow from a straightforward application of finite
geometric series
and the triangle inequality. Furthermore we also see that $x\in
\mbox{adh}_{\tau(s)}{\mathcal E}\cap [x_a-2^{-a+1},x_a].$

   We
next note that in this construction it is impossible that $x=x_a.$
Indeed otherwise $x_n=x_a$ whenever $n\in \mathbb{N}$ and $n\geq
a,$ and thus $x_a=x\in  \bigcap_{n\in \mathbb{N}} F_n\subseteq C,$
but we have chosen $x_a\not \in C.$ Hence we have proved our
claim.

We shall denote the Euclidean topology by $\tau(e)$ on $\mathbb{R}.$
   Choose $b\in \mathbb{R}$ such that $x<b<x_a.$
   Set $E:= \mbox{adh}_{\tau(e)}{\mathcal E}\cap [x_a-2^{-a+1},b].$
    Note that this set is nonempty, since $x$  belongs to it. We shall
show that any point
    $y\in E$
    is an accumulation point of $E$ with respect to the topology $\tau(e).$
Let $y\in E.$ In order to reach a contradiction suppose that
there is $q\in \mathbb{N}$ such that $(S^{-1}_{2^{-q}}(y)\cup
S_{2^{-q}}(y))\cap E=\{y\}$
where
without loss of generality we can assume that
   $(S^{-1}_{2^{-q}}(y)\cup S_{2^{-q}}(y))\subseteq
S^{-1}_{2^{-a+2}}(x_a).$
Because $y\not \in C,$ there is $m\in \mathbb{N}$ such that $y\not \in E_m.$
Since $y\in \mbox{adh}_{\tau(e)}{\mathcal E},$ we find $f_n\in E_n\cap
E_m \cap
(S^{-1}_{2^{-n}}(y)
\cup S_{2^{-n}}(y))$
   whenever $n\in \mathbb{N}$.
    In particular $f_n\not = y$ whenever $n\in \mathbb{N}.$

In the following we assume that $f_n\in S_{2^{-n}}(y)$ for infinitely
many $n\in \mathbb{N}$.
Denote  this subset of $\mathbb{N}$ by $L.$
(Otherwise we have that $f_n\in S^{-1}_{2^{-n}}(y)$ for infinitely many
$n\in \mathbb{N},$ a case
which can be treated
analogously by a conjugate method.)
Fix now $n\in \mathbb{N}$ such that $n > q.$
Furthermore consider any $p\in L$ with $p> n.$ Since ${\mathcal E}$
does not have a $\tau(e)$-cluster point in
$\mbox{cl}_{\tau(e)}S_{2^{-n}}(f_p)$, which is a subset of
$S_{2^{-q}}(y),$ we find $s_p\in \mathbb{N}$
such that $(\mbox{cl}_{\tau({e})}S_{2^{-n}}(f_p)) \cap
(\mbox{cl}_{\tau({e})} E_{s_p})=\emptyset$
by compactness
   of $\mbox{cl}_{\tau(e)}S_{2^{-n}}(f_p)$
    in the Euclidean topology $\tau(e)$ on $\mathbb{R}.$
Then indeed $S_{2^{{-n}}}(f_p)\cap F_{s_p}=\emptyset$ whenever $p\in
L$ and $p\geq n+1,$
   since $S_{2^{{-n}}}(f_p)$  is disjoint from $S_{2^{-a+1}}(C).$
Consequently $S_{2^{-(n+1)}}(f_p) \cap
S_{2^{-(n+1)}}^{-1}(F_{s_p})=\emptyset$ whenever
$p\in L$ and $p\geq n+1.$

Next we use a crucial general fact about the Sorgenfrey line: {\em
If $n\in \mathbb{N},$ $y\in \mathbb{R}$ and the sequence
$(b_m)_{m\in \mathbb{N}}$ converges to $y$ with respect to the
topology $\tau(e)$, then $]y,y+2^{-n}[ \subseteq \bigcup_{m\in
\mathbb{N}}[b_m,b_m+2^{-n}[$ and, analogously,
$]y-2^{-n},y[\subseteq \bigcup_{m\in
\mathbb{N}}]b_m-2^{-n},b_m].$}

Indeed let $z\in S_{2^{-n}}(y)$ and $z\not =y.$ Then $y\in S^{-1}_{2^{-n}}(z)$
   and so by an important property of the Sorgenfrey line, there is
$r\in \mathbb{N}$ such that
   $S^{-1}_{2^{-r}}(y)\cup S_{2^{-r}}(y)\subseteq S^{-1}_{2^{-n}}(z).$
   It follows that there is $m\in \mathbb{N}$ such that $b_m\in
S^{-1}_{2^{-r}}(y)\cup S_{2^{-r}}(y).$
   Thus $z\in S_{2^{-n}}(b_m)$. (The second part of the statement
   is established analogously.)

Therefore applying this argument to the sequence $(f_p)_{p\in L}$
   and its $\tau(e)$-limit $y,$
we see by the relationship established above that
$$]y,y+2^{-(n+1)}[\,\cap \bigcap_{p\in L,p\geq n+1}
   S_{2^{-(n+1)}}^{-1}(F_{s_p})=\emptyset.$$
By the definition of $F_{n+1}$ we then have that $F_{n+1}\subseteq
\bigcap_{h\in \mathbb{N}}S_{2^{-(n+1)}}^{-1}(F_h),$ which provides a contradiction,
since $f_{s}\in F_s\cap S_{2^{-s}}(y)$ for the infinitely many $s\in L,$ where
each $f_s\not= y.$
Hence each point in $E$
is an accumulation point of $E$ with respect to the topology $\tau(e).$

We now consider the nonempty closed subspace $E$ of the complete
metrizable space $(\mathbb{R},\tau(e))$.
Observe that ${\mathcal E}$ does not have a double
cluster point belonging to $E:$
   Otherwise, since ${\mathcal E}$ is finer than the filter
generated by $\{F_n:n\in \mathbb{N}\}$ on $\mathbb{R},$ this double cluster point must also be
a double cluster
   point of ${\mathcal F}.$ Hence it belongs to $C$, but we know that
the interval
$[x_a-2^{-a+1},b]$ is disjoint from $C.$

For each $n,m\in \mathbb{N}$ then set $A_{n,m}=\{x\in E: S^{-1}_{2^{-n}}(x)\cap
E_m=\emptyset\}.$
   Furthermore for each $n,m\in \mathbb{N}$ set $B_{n,m}=\{x\in E: S_{2^{-n}}(x)\cap
    E_m=\emptyset \}.$
   By our observation stated in the preceding paragraph $\{A_{n,m}:n,m\in \mathbb{N}\}
   \cup \{B_{n,m}:n,m\in \mathbb{N}\}$ is a cover of $E.$
   Hence by the Baire Category Theorem
   \cite[Theorem 3.9.3]{Engelking} applied to the subspace $E$ of
$(\mathbb{R},\tau(e)),$ we
   find a nonempty open
   real interval $I$ such that $\emptyset \not =(I\cap E)\subseteq
(\mbox{cl}_{\tau(e)}
    C_{n,m})\cap E$ for
   some $n,m\in \mathbb{N}$ where $C_{n,m}=A_{n,m}$
or $C_{n,m}=B_{n,m}.$ Let us consider in detail the second case.
The omitted argument for the first case is analogous.

   In the second case  we can conclude that any point $y$ belonging to
$I\cap E$ has
   a sequence $(b_k)_{k \in \mathbb{N}}$ in $E$ converging in $(\mathbb{R},\tau(e))$ to it such
that for each $k\in \mathbb{N},$ $S_{2^{-n}}(b_k)\cap
   E_m=\emptyset.$ By the crucial property of the Sorgenfrey line
discussed above,
   it follows  that
    $]y,y+2^{-n}[\,\cap
   E_m =\emptyset.$
   Choose $y_1\in I\cap E.$
Since the points of $E$ are not $\tau(e)$-isolated in $E$, we can find
a point $y_2\in I\cap E\cap (S^{-1}_{2^{-n}}(y_1)\cup
S_{2^{-n}}(y_1))$ distinct from
$y_1.$
Let $\alpha$ be the minimum of $\{y_1,y_2\}$ and $\beta$ be the maximum of
$\{y_1,y_2\}.$
We conclude that $\beta\in \, ]\alpha,\alpha+2^{-n}[$ and the latter set is disjoint
from $E_m.$
Hence $\beta\in E$ cannot be a $\tau(e)$-cluster point of ${\mathcal E}$,
which yields
   another contradiction. Hence we finally deduce that there is $n\in \mathbb{N}$ such
that $F_n\subseteq S_{2^{-a+3}}(C).$

Similarly it can be  shown that given $a\in \mathbb{N},$ the assumption that $F_n\setminus
S_{2^{-a+3}}^{-1}(C)\not=\emptyset$ whenever $n\in \mathbb{N}$ leads to a
contradiction.
We conclude by the K\"unzi-Ryser condition
 \cite[Proposition 8]{KunziRyser} that $({\mathcal
P}_0(\mathbb{R}),{({\mathcal U}_{s})}_H)$ is a bicomplete
quasi-uniform space.
 $\hfill{\Box}$

\section{The $2$-envelope of a filter}

The concepts of the envelope of a filter(base) and of a round filter on a quasi-uniform
space are well known (see e.g. \cite[p. 314]{Kunzi}). Similarly, in our context it is useful
to introduce the concepts of a $2$-envelope of a filter and of a $2$-round filter.

\begin{Def}
Let $(X,{\mathcal U})$ be a quasi-uniform space and let ${\mathcal F}$ be
a filterbase on $X.$
Then we consider the filter ${\mathcal D}_{{\mathcal U}}({\mathcal F})$
generated by the base $\{U^{-1}(F)\cap U(F):F\in {\mathcal F},U\in
   {\mathcal U}\}$ on $X.$
It will be called the {\em $2$-envelope} of ${\mathcal F}.$
A filter is called {\em $2$-round} if it is equal to its $2$-envelope.
\end{Def}

\begin{Lem} Let ${\mathcal F}$ be a doubly stable filter on a
quasi-uniform space $(X,{\mathcal U}).$
Then ${\mathcal D}_{\mathcal U}({\mathcal F})$ is doubly stable.
\end{Lem}

{\em Proof.} Let $U\in {\mathcal U}.$ Choose $V\in {\mathcal U}$ such that
   $V^2\subseteq U.$
There is $F_V\in {\mathcal F}$ such that for each $F\in {\mathcal F}$ we have
$F_V\subseteq V^{-1}(F)\cap V(F).$ So for all $F\in {\mathcal F}$ and all
   $W\in {\mathcal U}$
   we get $F_V\subseteq V^{-1}(W^{-1}(F)\cap W(F))\cap V(W^{-1}(F)\cap W(F)).$
Therefore for each $F\in {\mathcal F}$ and $W\in {\mathcal U}$ we have
$V^{-1}(F_V)\subseteq V^{-2}(W^{-1}(F)\cap W(F))$ and $V(F_V)\subseteq
V^{2}(W^{-1}(F)\cap W(F)).$
Hence for all $F\in {\mathcal F}$ and $W\in {\mathcal U}$, we see that
$V^{-1}(F_V)\cap V(F_V)\subseteq U(W^{-1}(F)\cap W(F))$ and
$V^{-1}(F_V)\cap V(F_V)\subseteq U^{-1}(W^{-1}(F)\cap W(F)).$
We have shown that ${\mathcal D}_{\mathcal U}({\mathcal F})$ is
doubly stable.
$\hfill{\Box}$

\begin{Lem}\label{zwei}  For any filter ${\mathcal F}$ on a quasi-uniform space
$(X,{\mathcal U})$,
   ${\mathcal D}_{\mathcal U}({\mathcal F})$ is
2-round, that is, we
have ${\mathcal D}_{\mathcal U}({\mathcal D}_{\mathcal U}({\mathcal F}))=
{\mathcal D}_{\mathcal U}({\mathcal F}).$
Furthermore ${\mathcal D}_{\mathcal U}({\mathcal F})$
has the same sets of
$\tau({\mathcal U})$-
and $\tau({\mathcal U}^{-1})$-cluster points as ${\mathcal F}.$
\end{Lem}

{\em Proof.} Clearly
${\mathcal D}_{\mathcal U}({\mathcal F})\subseteq {\mathcal F}.$
Therefore ${\mathcal D}_{\mathcal U}({\mathcal D}_{\mathcal U}
({\mathcal F}))\subseteq {\mathcal D}_{\mathcal U} ({\mathcal F}).$
Let $U\in {\mathcal U}.$ Choose $V\in {\mathcal U}$ such that
$V^2\subseteq U.$
Then for any $F\in {\mathcal F}$ we have $V^{-1}(V^{-1}(F)\cap
V(F))\cap V(V^{-1}(F)\cap V(F))
\subseteq V^{-2}(F)\cap V^2(F)\subseteq U^{-1}(F)\cap U(F).$
Thus ${\mathcal D}_{\mathcal U}({\mathcal F})\subseteq
{\mathcal D}_{\mathcal U}({\mathcal D}_{\mathcal U}({\mathcal F}))$ and
so the assertion holds.

Since ${\mathcal D}_{\mathcal U}({\mathcal F})$ is coarser than ${\mathcal F},$
the filter ${\mathcal D}_{\mathcal U}({\mathcal F})$
 certainly has all the $\tau({\mathcal U})$- and $\tau({\mathcal U}^{-1})$-cluster
 points of ${\mathcal F}.$
On the other hand it is evident by the definition of the generating
filterbase of
${\mathcal D}_{\mathcal U}({\mathcal F})$ that if $x\in X$ is, say,
   a $\tau({\mathcal U})$-cluster point
of ${\mathcal D}_{\mathcal U}({\mathcal F}),$ then $x$ is also a
   $\tau({\mathcal U})$-cluster point of ${\mathcal F}:$
   Indeed let $U\in {\mathcal U}.$ Then $U(x)\cap U^{-1}(F)\not=\emptyset$
whenever $F\in {\mathcal F}$
implies that $U^2(x)\cap F\not =\emptyset$
 whenever $F\in {\mathcal F}.$ The corresponding result obviously also holds for
 $\tau({\mathcal U}^{-1})$-cluster points.
 $\hfill{\Box}$

\begin{Rem}
Given a quasi-uniform space $(X,{\mathcal U})$, note that for any $C\in {\mathcal P}_0(X),$
we obviously have   ${\mathcal D}_{\mathcal U}({\mathcal C}_C)=({\mathcal
C}_{C})_{\mathcal U}.$
 \footnote{We are going to show that this equality holds for an arbitrary doubly stable filter on
a quasi-uniform space.}
Furthermore $({\mathcal C}_C)_{\mathcal U}=
({\mathcal C}_{C'})_{\mathcal U}$ where $C'$ is the double closure
of $C$ in $(X,{\mathcal U}).$
$\hfill{\Box}$
\end{Rem}

\begin{Lem} \label{drei} Let ${\mathcal F}$ be a doubly stable filter
on a quasi-uniform
space $(X,{\mathcal U}).$ Then ${\mathcal F}_{\mathcal U}=
{\mathcal D}_{\mathcal U}({\mathcal F}).$
\end{Lem}

{\em Proof.} By definition, ${\mathcal F}_{\mathcal U}$ is
generated by the base $\{\bigcap_{F\in {\mathcal F}}(U^{-1}(F)\cap
U(F)):U\in {\mathcal U}\}$, while ${\mathcal D}_{\mathcal
U}({\mathcal F})$ is generated by $\{U^{-1}(F)\cap U(F):U \in
{\mathcal U},F\in {\mathcal F}\}.$ Thus clearly ${\mathcal
D}_{\mathcal U}({\mathcal F})\subseteq {\mathcal F}_{\mathcal U}.$
In order to establish equality, let $W\in {\mathcal U}$ and choose
$U\in {\mathcal U}$ such that $U^3\subseteq W.$ Then $U_{\mathcal
F}\subseteq \bigcap_{F\in {\mathcal F},V\in {\mathcal U}}
U^{-1}(V^{-1}(F)\cap V(F))\cap \bigcap_{F\in {\mathcal F},V\in
{\mathcal U}} U(V^{-1}(F)\cap V(F)).$ Thus $U^{-1}(U_{\mathcal
F})\cap U(U_{\mathcal F})\subseteq $ $$U^{-1} (\bigcap_{F\in
{\mathcal F},V\in {\mathcal U}} U^{-1} (V^{-1}(F)\cap V(F)))\cap
U(\bigcap_{F\in {\mathcal F},V\in {\mathcal U}}U (V^{-1}(F)\cap
V(F)))\subseteq $$ $\bigcap_{F\in {\mathcal F}}(U^{-3}(F)\cap
U^3(F)) \subseteq \bigcap_{F\in {\mathcal F}}(W^{-1}(F)\cap
W(F))=W_{\mathcal F}.$ By the last chain of inclusions we conclude
that ${\mathcal F}_{\mathcal U}\subseteq {\mathcal D}_{\mathcal
U}({\mathcal F}),$ since ${\mathcal F}$ is doubly stable and so ${
U}_{\mathcal F}\in {\mathcal F}.$ $\hfill{\Box}$

\begin{Cor} \label{base} For each doubly stable filter ${\mathcal F}$
on a quasi-uniform
space $(X,{\mathcal U})$ the filter ${\mathcal F}_{\mathcal U}$ has a base consisting
of $\tau({\mathcal U}^s)$-open subsets of $X.$
\end{Cor}

{\em Proof.} Obviously $\{\mbox{int}_{\tau({\mathcal U}^{-1})}U^{-1}(F)\cap
\mbox{int}_{\tau({\mathcal U})}U(F):U\in {\mathcal U},F\in {\mathcal F}\}$
is such a base for ${\mathcal D}_{\mathcal U}({\mathcal F}),$ since for
instance for some $U,W\in {\mathcal U}$ with $W^2\subseteq U$ we
have $W^{-1}(F)\subseteq \mbox{int}_{\tau({\mathcal U}^{-1})}U^{-1}(F).$
The assertion now follows from Lemma \ref{drei}. $\hfill{\Box}$

\section{The main construction}

In this section
 we introduce a stability functor on the category of quasi-uniform
spaces and quasi-uniformly continuous maps and compare it with
 the Hausdorff hyperspace
functor and the bicompletion functor.
The definition contained in our next proposition is obviously motivated
by the construction of the Hausdorff quasi-uniformity (see e.g. \cite{KunziRyser}).

\begin{Prp} \label{stabcon}
Let $(X,{\mathcal U})$ be a quasi-uniform space and let $S_D(X)$ be the set of all
doubly stable filters on $(X,{\mathcal U}).$

For each $U\in {\mathcal U}$ we set

${U}_+=\{({\mathcal F},{\mathcal G})\in S_D(X)\times S_D(X):
\bigcap_{F\in {\mathcal F}}U(F)\in {\mathcal G}\}.$

Then $\{{U}_+: U\in {\mathcal U}\}$ is a base for the {\em
upper quasi-uniformity}
${\mathcal U}_+$
   on $S_D(X).$

For each $U\in {\mathcal U}$ we set

${U}_-=\{({\mathcal F},{\mathcal G})\in S_D(X)\times S_D(X):$
   $\bigcap_{G\in {\mathcal G}} U^{-1}(G)\in {\mathcal F}\}.$

Then $\{U_-: U\in {\mathcal U}\}$
   is a base for the {\em lower quasi-uniformity} ${\mathcal U}_-$
   on $S_D(X).$

Furthermore for each $U\in {\mathcal U}$ set $U_{D}=U_+\cap U_-.$
Then $\{U_{D}:U\in {\mathcal U}\}$ is a base for the
{\em stability
quasi-uniformity} ${\mathcal U}_{D}$
on $S_D(X).$
The {\em stability space} of $(X,{\mathcal U})$ is the $T_0$-quotient space of
$(S_D(X),{\mathcal U}_D)$ and will be denoted by $(qS_D(X),q{\mathcal U}_D).$

\end{Prp}

{\em Proof.}
Note first that for each $U\in {\mathcal U}$ and any
 ${\mathcal F}\in S_D(X),$
we have $({\mathcal
F},{\mathcal F})\in
{U}_+,$ and similarly $({\mathcal F},{\mathcal F})\in U_-$ and
$({\mathcal F},{\mathcal F})\in U_D.$
Observe also that $U,V\in {\mathcal U}$ with $U\subseteq V$ implies that
$U_+\subseteq V_+,$ $U_-\subseteq V_-$, and ${U}_D\subseteq {V}_D.$
Hence $\{{U}_+:U\in {\mathcal U}\}$, $\{{U}_-:U\in {\mathcal U}\}$ and
$\{{U}_D:U\in {\mathcal U}\}$  are filterbases on $S_D(X)\times S_D(X).$

Let $U\in {\mathcal U}$ and $V\in {\mathcal U}$ be such that
$V^{2}\subseteq U.$
Let $({\mathcal F},{\mathcal G})\in {V}_+$
and $({\mathcal G},{\mathcal H})\in {V}_+.$
Then there is $G\in {\mathcal G}$ such that $G\subseteq V(F)$
whenever $F\in {\mathcal F}.$
Similarly there is $H\in {\mathcal H}$ such that $H\subseteq V(G)$
whenever $G\in {\mathcal G}.$
Hence there is $H\in {\mathcal H}$ such that $H\subseteq U(F)$
whenever $F\in {\mathcal F}.$
We have shown that $({\mathcal F},{\mathcal H})\in {U}_+.$
Thus $(V_+)^2\subseteq U_+.$ Similarly $(V_-)^2\subseteq U_-,$ and thus
$(V_D)^2\subseteq U_D.$

We deduce that ${\mathcal U}_+,$ ${\mathcal U}_-,$ and ${\mathcal U}_D$ are
   quasi-uniformities on $S_D(X).$ $\hfill{\Box}$

 \begin{Rem} Let $(X,{\mathcal U})$ be a quasi-uniform space.
 Then $(S_D(X),({\mathcal U}_D)^{-1})=(S_D(X),({\mathcal U}^{-1})_D).$
\end{Rem}

{\em Proof.} For any $U\in {\mathcal U}$ we have $(U_+)^{-1}=(U^{-1})_-$
and $(U_-)^{-1}=(U^{-1})_+$ and thus $(U_D)^{-1}=(U^{-1})_D.$
The assertion follows. $\hfill{\Box}$

\begin{Cor} Let $(X,{\mathcal U})$ be a quasi-uniform space and let $S_D(X)$
be the set of all doubly stable filters on $X.$
Then $({\mathcal U}_+)^{-1}=({\mathcal U}^{-1})_-$ and $({\mathcal U}_-)^{-1}=
({\mathcal U}^{-1})_+$ on $S_D(X).$ $\hfill{\Box}$

\end{Cor}

\begin{Prp} Let $(X,{\mathcal U})$ and $(Y,{\mathcal V})$ be quasi-uniform
spaces and $f:(X,{\mathcal U})\rightarrow (Y,{\mathcal V})$ be a
quasi-uniformly
continuous map.

If ${\mathcal F}$ is a doubly stable filter on $(X,{\mathcal U})$,
then the filter $[f({\mathcal F})]$ generated by the filterbase
$\{f(F):F\in {\mathcal F}\}$ on $Y$
is doubly stable on $(Y,{\mathcal V}).$
Furthermore the map
${f}_D:(S_D(X),{\mathcal U}_D)\rightarrow (S_D(Y),{\mathcal V}_D)$
defined by ${f}_D({\mathcal F})=[f({\mathcal F})]$ is quasi-uniformly
continuous.
\end{Prp}

{\em Proof.} Let $V\in {\mathcal V}.$ By quasi-uniform continuity of $f$
there is $U\in {\mathcal U}$ such that
$(f\times f)U\subseteq V.$
Since ${\mathcal F}$ is doubly stable,
there is $F_U\in {\mathcal F}$ such that
$F_U\subseteq U^{-1}(F)\cap U(F)$ whenever $F\in {\mathcal F}.$ Consequently
$f(F_U)\subseteq V^{-1}(f(F))\cap V(f(F))$
whenever $F\in {\mathcal F}.$
We conclude that $[f({\mathcal F})]$ is
doubly stable on $(Y,{\mathcal V})$ and
thus ${f}_D$ is well-defined.

It remains to show that ${f}_D$ is quasi-uniformly continuous.
Let $V\in {\mathcal V}.$ As above,
there is $U\in {\mathcal U}$ such that
$(f\times f)U\subseteq V.$
Consider $({\mathcal F},{\mathcal G})\in {U}_D.$
Then $\bigcap_{F\in {\mathcal F}}U(F)\in {\mathcal G}$ and
$\bigcap_{G\in {\mathcal G}}
U^{-1}(G)\in{\mathcal F}.$
Consequently $\bigcap_{F\in {\mathcal F}}V(f(F))\in
[f({\mathcal G})]$ and
$\bigcap_{G\in {\mathcal G}}V^{-1}(f(G))\in [f({\mathcal F})],$ and thus
$({f}_D({\mathcal F}),{f}_D({\mathcal G}))\in {V}_D.$
Hence the map ${f}_D$
is quasi-uniformly continuous. $\hfill{\Box}$

\begin{Rem} Given a quasi-uniform space $(X,{\mathcal U})$
several authors (see e.g.  \cite{Deak,Doitchinovstable}) have
considered kinds of extensions of $(X,{\mathcal U})$ based on the
concept of envelopes in $(X,{\mathcal U}).$ Often these
constructions can be understood as generalizations of our
construction in Proposition \ref{stabcon} above to more general
collections of filters or  even families of filter pairs.

In fact, given for instance
 any collection ${\mathcal M}$ of (round) filters on a quasi-uniform space $(X,{\mathcal U})$ a
quasi-uniformity ${\mathcal U}_\oplus$
on ${\mathcal M}$ can be defined which has $\{U_\oplus:U\in {\mathcal U}\}$
as a base, where ${U}_\oplus=\{({\mathcal F},{\mathcal G})\in {\mathcal M}\times {\mathcal M}:
U(F)\in {\mathcal G}$ whenever $F\in {\mathcal F}\}.$
Similarly we can define on ${\mathcal M}$ a quasi-uniformity ${\mathcal U}_\ominus$
generated by the base $\{U_\ominus:U\in {\mathcal U}\}$
 where ${U}_\ominus:=(({ U}^{-1})_\oplus)^{-1}$ whenever $U\in {\mathcal U}.$

In case that ${\mathcal M}$ only consists of stable filters, one readily verifies that
for each $U\in{\mathcal U}$ we have $U_+\subseteq {U}_\oplus$ and  ${U}_\oplus\subseteq
(U^2)_+.$
Hence indeed ${\mathcal U}_+={{\mathcal U}_\oplus}$ on ${\mathcal M}.$ Additionally
we also have
${\mathcal U}_-={\mathcal U}_\ominus$ if ${\mathcal M}$ even consists of doubly stable filters.
$\hfill{\Box}$
\end{Rem}

\begin{Rem} Given a quasi-uniform space $(X,{\mathcal U}),$ we remark
 that if $({\mathcal F}, {\mathcal G})$ is a Cauchy filter pair on $(X,{\mathcal U})$ in the
sense of Doitchinov\footnote{A pair $({\mathcal F},{\mathcal G})$
 of filters on a quasi-uniform space $(X,{\mathcal U})$ is called {\em a Cauchy filter pair}
  if for each $U\in {\mathcal U}$
there are $F\in {\mathcal F}$ and $G\in {\mathcal G}$
such that $F\times G\subseteq U$
\cite{Doitchinov}.}, then
$({\mathcal F},{\mathcal G})\in (\bigcap {{\mathcal U}}_\oplus) \cap (\bigcap {\mathcal U}_\ominus).$
\end{Rem}

{\em Proof.} Let $U\in {\mathcal U}.$ Then there are $F\in {\mathcal
F}$ and $G\in {\mathcal G}$
such that $F\times G\subseteq U.$ Thus $G\subseteq U(F')$ whenever
$F'\in {\mathcal F},$ since
$F'\cap F\not =\emptyset.$
Similarly $F\subseteq U^{-1}(G')$ whenever $G'\in {\mathcal G}$, since
$G'\cap G\not = \emptyset.$
Thus $({\mathcal F},{\mathcal G})\in U_\oplus\cap U_\ominus.$
$\hfill{\Box}$

\bigskip
The preceding remark suggests to call
a filter pair $( {\mathcal F},{\mathcal G})$
of doubly
stable filters  on
a quasi-uniform space $(X,{\mathcal U})$
{\em generalized Cauchy} provided that $({\mathcal F},{\mathcal G})\in
\bigcap {\mathcal U}_D.$ In this article however there will be no need
to study this concept further.

\begin{Rem}
Let us note that on the subset $\widetilde{X}$ of $S_D(X)$
consisting of all the minimal ${\mathcal U}^s$-Cauchy filters,
our definition of ${\mathcal U}_+$ (resp. ${\mathcal U}_-$)
 yields the standard (explicit) construction of the bicompletion quasi-uniformity $\widetilde{\mathcal U}$
 \cite[Theorem 3.33]{FletcherLindgren} of a quasi-uniform $T_0$-space $(X,{\mathcal U})$,\footnote{The quasi-uniformity
$\widetilde{\mathcal U}$ on $\widetilde{X}$ is generated by the base $\{\widetilde{U}:U\in {\mathcal U}\}$
where for any $U\in {\mathcal U}$
we have $\widetilde{U}=\{({\mathcal F},{\mathcal G})\in \widetilde{X}\times \widetilde{X}:$
there are $F\in {\mathcal F}$
and $G\in {\mathcal G}$ such that $F\times G\subseteq U\}.$} as we show next:

Fix $U\in {\mathcal U}.$
For minimal ${\mathcal U}^s$-Cauchy filters ${\mathcal F}$ and ${\mathcal G}$
on $(X,{\mathcal U})$ suppose that there are $F\in {\mathcal F}$
   and $G\in {\mathcal G}$ such that $F\times G\subseteq U.$
By the argument given above,
 $\bigcap_{F'\in {\mathcal F}}U(F') \in {\mathcal G}$
and $\bigcap_{G'\in {\mathcal G}}U^{-1}(G')\in {\mathcal F}$ so that
$({\mathcal F},{\mathcal G})\in U_+\cap U_-.$

On the other hand suppose that
$V\in {\mathcal U}$ such that $V^2\subseteq U.$
Let $\bigcap_{F\in {\mathcal F}}V(F)\in
{\mathcal G}$
(resp. $\bigcap_{G\in {\mathcal G}}V^{-1}(G)\in {\mathcal F}).$
   Furthermore since ${\mathcal F}$ and ${\mathcal G}$
    are (minimal) ${\mathcal U}^s$-Cauchy filters, there exist $F'\in {\mathcal
F}$ and $G'\in {\mathcal G}$ such that
$(F'\times F')\cup (G'\times G')\subseteq V.$
It follows that $F'\times V(F')\subseteq V^2$ and $V^{-1}(G')\times G'\subseteq V^2.$

Consequently
$F'\times \bigcap_{F\in {\mathcal F}}V(F)\subseteq U$
(resp. $\bigcap_{G\in {\mathcal G}}V^{-1}(G)\times G'\subseteq U).$
Therefore in either case  there are $F\in {\mathcal F}$ and $G\in {\mathcal G}$
such that $F\times G\subseteq U$ and the
claim is verified.

By the same argument we conclude that for any quasi-uniform $T_0$-space
$(X,{\mathcal U})$ the subspace $(\widetilde{X},{\mathcal U}_D\vert \widetilde{X})$
of $(S_D(X),{\mathcal U}_D)$ is quasi-uniformly isomorphic to
 the bicompletion $(\widetilde{X},\widetilde{\mathcal U})$ of $(X,{\mathcal U})$
with the quasi-uniform embedding defined by $x\mapsto {\mathcal D}_{\mathcal U}
({\mathcal C}_{\{x\}})= {\mathcal U}^s(x)$ whenever
$x\in X.$
\end{Rem}

\begin{Rem} Let $(X,{\mathcal U})$ be a quasi-uniform space.

Suppose that
$S_D(X^s)$ denotes the set of all ${\mathcal U}^s$-stable filters on $X.$
Of course, $S_D(X^s)\subseteq S_D(X),$ and $S_D(X^s)$ is the carrier set of the
uniformity $({\mathcal U}^s)_D.$

Furthermore $({\mathcal U}_D)^s\vert S_D(X^s)\subseteq ({\mathcal U}^s)_D.$
\end{Rem}

{\em Proof.} Indeed $({\mathcal U}_D)^s$ restricted to $S_D(X^s)$ is generated
by the base consisting of all entourages
$(U_D)^s \cap (S_D(X^s)\times S_D(X^s))=\{({\mathcal F},{\mathcal G})\in S_D(X^s)
\times S_D(X^s):
U_{\mathcal F}\in {\mathcal G}$ and
$U_{\mathcal G}\in {\mathcal F}\}$ where $U\in {\mathcal U},$
while the set of all entourages
$(U^s)_D=\{({\mathcal F},{\mathcal G})\in S_D(X^s)
\times S_D(X^s):
(U^s)_{\mathcal F}\in {\mathcal G}$ and $(U^s)_{\mathcal G}\in {\mathcal F}\}$
with $U\in {\mathcal U}$
generates $({\mathcal U}^s)_D.$ The assertion follows.
$\hfill{\Box}$

\begin{Rem} \label{categ}
(a) The map ${\mathcal C}(C)= {\mathcal C}_{C}$
for any $C\in {\mathcal P}_0(X)$ defines a quasi-uniform embedding of the
Hausdorff hyperspace
$({\mathcal P}_0(X),{\mathcal U}_H)$
into the quasi-uniform space $(S_D(X),{\mathcal U}_D).$

(b) For any quasi-uniformly continuous map $f:(X,{\mathcal
U})\rightarrow
(Y,{\mathcal V})$ between quasi-uniform spaces $(X,{\mathcal U})$ and
$(Y,{\mathcal V})$, the map
$f_D:(S_D(X),{\mathcal U}_D)\rightarrow (S_D(Y),{\mathcal V}_D)$
restricts to the usual hypermap ${\mathcal
P}_0(f):({\mathcal
P}_0(X),{\mathcal U}_H)\rightarrow ({\mathcal P}_0(Y),{\mathcal V}_H)$
where according to part (a) the two hyperspaces are considered as the subspaces
${\mathcal C}{\mathcal P}_0(X)$ and ${\mathcal C}{\mathcal P}_0(Y)$
of $(S_D(X),{\mathcal U}_D)$ resp. $(S_D(Y),{\mathcal V}_D).$
\end{Rem}

   {\em Proof.} (a)
   Clearly ${\mathcal C}$ is injective. We verify that it is a quasi-uniform
   embedding:
    For any $U\in {\mathcal U},$ we have that $(A,B)\in
U_H$ if, and only if $B\subseteq
   U(A)$ and $A\subseteq U^{-1}(B)$
if, and only if $({\mathcal C}_A,{\mathcal C}_B)\in {U}_D.$

(b) Of course the usual hypermap on ${\mathcal P}_0(X)$ into ${\mathcal P}_0(Y)$
is defined by $A\mapsto f(A)$.
Indeed the restriction to
${\mathcal C}{\mathcal P}_0(X)$ of our map $f_D$ is given by
${\mathcal C}_A\mapsto {\mathcal C}_{fA}.$
$\hfill{\Box}$

\begin{Lem}  \label{rep} Two doubly stable filters ${\mathcal F}$ and ${\mathcal G}$ on
a quasi-uniform space $(X,{\mathcal U})$ are ${\mathcal U}_D$-equivalent
if and only if ${\mathcal F}_{\mathcal U}={\mathcal G}_{\mathcal U}.$
(Hence for each
   doubly stable filter ${\mathcal F}$ on a quasi-uniform space $(X,{\mathcal U}),$
   the doubly stable filter ${\mathcal F}_{\mathcal U}$
can represent its ${{\mathcal U}}_D$-equivalence class on $S_D(X)$.)
\end{Lem}

   {\em Proof.} Indeed if ${\mathcal F}$ and ${\mathcal G}$
    are ${\mathcal U}_D$-equivalent, then by definition, ${\mathcal F}_{\mathcal
U}\subseteq {\mathcal G}$
     and ${\mathcal G}_{\mathcal U}\subseteq {\mathcal F}.$
So by monotonicity of the $_{\mathcal U}$-operator, $({\mathcal
F}_{\mathcal U})_{\mathcal U}\subseteq {\mathcal G}_{\mathcal U}$
and $({\mathcal G}_{\mathcal U})_{\mathcal U}\subseteq {\mathcal
F}_{\mathcal U}.$
Therefore ${\mathcal F}_{\mathcal U}={\mathcal G}_{\mathcal U},$
because the operation
$\cdot_{\mathcal U}$ is idempotent by Lemmas \ref{zwei} and \ref{drei}.

If ${\mathcal F}_{\mathcal U}=
{\mathcal G}_{\mathcal U},$ then $U_{\mathcal F}\in {\mathcal G}$
and $U_{\mathcal G}\in {\mathcal F}$ whenever $U\in {\mathcal U}$.
So ${\mathcal F}$ and ${\mathcal G}$ are ${{\mathcal U}_D}$-equivalent by
the definition of this equivalence relation.
$\hfill{\Box}$

\begin{Lem}\label{dense} Let $(X,{\mathcal U})$ be a quasi-uniform
space and let $A$ be a
$\tau({\mathcal U}^s)$-dense subspace of $X.$ If ${\mathcal F}$ is a
doubly stable filter on
$(X,{\mathcal U}),$
then $\{U_{\mathcal F}\cap A:U\in {\mathcal U}\}$
is a filterbase of a  doubly stable filter
${\mathcal F}_{\mathcal U}\vert A$ on $(A,{\mathcal U}\vert A).$

 The filter $[{\mathcal F}_{\mathcal U}\vert A]$ is
${\mathcal U}_D$-equivalent to ${\mathcal F},$ where
 $[{\mathcal F}_{\mathcal U}\vert A]$ denotes the filter generated on $X$
by the filterbase ${\mathcal F}_{\mathcal U}\vert A.$
\end{Lem}

{\em Proof.} We first mention that ${\mathcal F}_{\mathcal U}$ has
a base of  $\tau({\mathcal U}^s)$-open
sets, see Corollary \ref{base}.
So the filter ${\mathcal F}_{\mathcal U}\vert A$ is well defined.
Recall also \cite[Theorem 1.3.6]{Engelking} that if $A$ is dense and $G$ open in a topological space, then
$\overline{G\cap A}=\overline{G}.$
For any $U\in {\mathcal U},$ by definition we have
$U_{{\mathcal F}_{\mathcal U}}=\bigcap_{F\in {\mathcal F}_{\mathcal
U}}(U(F)\cap
U^{-1}(F)),$
while $U_{[{\mathcal F}_{\mathcal U}\vert A]}=\bigcap_{F\in
{\mathcal F}_{\mathcal U}}(U(F\cap A)\cap U^{-1}(F\cap A)).$
Therefore clearly
   $U_{[{\mathcal F}_{\mathcal U}\vert A]}\subseteq U_{{\mathcal
F}_{\mathcal U}}.$

Let $U,V\in {\mathcal U}$ be such that $V^2\subseteq U.$
We check that
   $V_{{\mathcal F}_{\mathcal U}}\subseteq U_{[{\mathcal F}_{\mathcal U}
\vert A]}:$
Indeed $$V_{{\mathcal F}_{\mathcal U}}\subseteq
\bigcap_{F\in {\mathcal
F}_{\mathcal U}\cap \tau({\mathcal U}^s)}
(V^{-1}(\mbox{cl}_{\tau({\mathcal U}^{s})}F)\cap
V(\mbox{cl}_{\tau({\mathcal U}^s)}{F}))\subseteq $$ $$
\bigcap_{F\in {\mathcal F}_{\mathcal U}\cap \tau({\mathcal U}^s)}
(V^{-1}(\mbox{cl}_{\tau({\mathcal U}^{s})}(F\cap A))\cap
V(\mbox{cl}_{\tau({\mathcal U}^s)}(F\cap A))\subseteq$$ $$
\bigcap_{F\in {\mathcal F}_{\mathcal U}}(V^{-2}(F\cap A)\cap
V^{2}(F\cap A))\subseteq
(V^{2})_{[{\mathcal F}_{{\mathcal U}}
\vert_A]}\subseteq U_{[{\mathcal F}_{{\mathcal U}\vert A}]},$$
where we have used that ${\mathcal F}_{\mathcal U}$ has a base of
 $\tau({\mathcal U}^s)$-open sets.

Since $V_{{{\mathcal F}_{\mathcal U}}}\in {\mathcal F}_{{\mathcal
U}}$ and so $V_{{{\mathcal F}_{\mathcal U}}} \cap A \in {\mathcal
F}_{\mathcal U}\vert A,$ by the last chain of inequalities we
first conclude that ${\mathcal F}_{\mathcal U}\vert A$ is doubly
stable on $(A,{\mathcal U}\vert A).$ Hence by Lemma \ref{subs}
$[{\mathcal F}_{\mathcal U}\vert A]$ is doubly stable on
$(X,{\mathcal U}).$ Furthermore by the argument just presented we
also see that the filters
   ${\mathcal F}_{\mathcal U}$ and
$[{\mathcal F}_{\mathcal U}
\vert A]$ are ${\mathcal U}_D$-equivalent,
that is,
${\mathcal F}_{\mathcal U}
=[{\mathcal F}_{\mathcal U}\vert A]_{\mathcal U}.$ $\hfill{\Box}$

 \begin{Cor} \label{resbic}
Let $(X,{\mathcal U})$ be a quasi-uniform
$T_0$-space and $(\widetilde{X},\widetilde{\mathcal U})$
its bicompletion with quasi-uniform embedding $i:(X,{\mathcal U})\rightarrow
(\widetilde{X},\widetilde{\mathcal U})$ where
$i(x)={\mathcal U}^s(x)$ whenever $x\in X.$ (To simplify the notation in the following
we shall often identify each $x\in X$ with $i(x)$ and consider $(X,{\mathcal U})$
as a subspace of $(\widetilde{X},\widetilde{\mathcal U}).$
)
If ${\mathcal F}$ is a doubly stable filter on
$(\widetilde{X},\widetilde{\mathcal U})$,
then $[{\mathcal F}_{\widetilde{\mathcal U}}\vert X]
_{\widetilde{\mathcal U}}={\mathcal F}
_{\widetilde{\mathcal U}}.$
\end{Cor}

{\em Proof.} The assertion follows from the preceding result. $\hfill{\Box}$

\bigskip
Let $(X,{\mathcal U})$ be a quasi-uniform $T_0$-space
and $(\widetilde{X},\widetilde{\mathcal U})$ its bicompletion.
We shall now consider the following commutative diagram,
where the maps ${\mathcal C}$, $\widetilde{{\mathcal C}},$ $e_0$ and $e$
are defined as follows:

First ${\mathcal C}(C)={\mathcal C}_C$ is equal to the filter generated by the base $\{C\}$
on $X$ whenever $C\in {\mathcal P}_0(X).$
Moreover $\widetilde{\mathcal C}(C)=\widetilde{{\mathcal C}}_C$ is equal to
the filter generated by the base $\{C\}$
on $\widetilde{X}$ whenever $C\in {\mathcal P}_0(\widetilde{X}).$
We remark that either map is the quasi-uniform embedding described in Remark \ref{categ}.

Furthermore $e_0(C)=C$ whenever $C\in {\mathcal P}_0(X).$
Finally $e({\mathcal F})$ is equal to the filter $[{\mathcal F}]$
generated on $\widetilde{X}$ by the filterbase ${\mathcal F},$ where ${\mathcal F}$ is
a doubly stable filter on $(X,{\mathcal U}).$
It is readily checked that $e_0$ and $e$ are quasi-uniformly continuous.
 Indeed we have $e=i_D$ for
the usual quasi-uniform embedding $i:(X,{\mathcal U})\rightarrow (\widetilde{X},
\widetilde{\mathcal U}).$\footnote{In the proof of
Proposition \ref{bc} we shall sketch a direct argument that establishes
quasi-uniform continuity of the related map $qe.$}

$$\begin{xy}
\xymatrix{
({\mathcal P}_0(X),{\mathcal U}_H) \ar[r]^{\mathcal C} \ar[d]_{e_0} & (S_D(X),{\mathcal U}_D)
 \ar[d]^{e} \\
({\mathcal P}_0(\widetilde{X}),\widetilde{\mathcal U}_H)
\ar[r]_{\widetilde{\mathcal C}} & (S_D(\widetilde{X}),\widetilde{\mathcal U}_D)
}
\end{xy}
$$

Applying the $T_0$-reflector to our first diagram yields another commutative diagram
for the corresponding $T_0$-quotient
spaces (see our second diagram below).
We shall interpret $q{\mathcal P}_0(X),$ resp. $q{\mathcal P}_0(\widetilde{X}),$
as the subspaces of $({\mathcal P}_0(X),{\mathcal U}_H),$ resp.
$({\mathcal P}_0(\widetilde{X}),\widetilde{\mathcal U}_H),$
consisting of all nonempty
 doubly closed subsets of $(X,{\mathcal U})$ resp. $(\widetilde{X},\widetilde{\mathcal U})$,
 and $qS_D(X)$ resp. $qS_D(\widetilde{X})$
as the subspaces of $(S_D(X),{\mathcal U}_D)$ resp.
$(S_D(\widetilde{X}),\widetilde{\mathcal U}_D)$ consisting of all
the filters ${\mathcal F}_{\mathcal U}$ resp. $({\mathcal
F}')_{\widetilde{\mathcal U}},$ where ${\mathcal F}$ is a doubly
stable filter on $(X,{\mathcal U})$ resp. ${\mathcal F}'$ is a
doubly stable filter on $(\widetilde{X},\widetilde{\mathcal U})$
(compare Theorem \ref{main} below).

Hence $q{\mathcal C}(C)$ is the filter on $X$ generated by the
base $\{U^{-1}(C)\cap U(C): U\in {\mathcal U}\}$ where $C$ is
doubly closed in $(X,{\mathcal U}),$ and $q\widetilde{{\mathcal
C}}(C')$ is the filter on $\widetilde{X}$ generated by the base
$\{\widetilde{U}^{-1}(C')\cap \widetilde{U}(C'):U\in {\mathcal
U}\}$ where $C'$ is doubly closed in
$(\widetilde{X},\widetilde{\mathcal U}).$ Furthermore $qe_0 (C)=
\mbox{cl}_{\tau(\widetilde{{\mathcal U}})}C\cap
\mbox{cl}_{\tau(\widetilde{{\mathcal U}}^{-1})}C$ where $C$ is
doubly closed in $(X,{\mathcal U}).$ Moreover $qe({\mathcal
F}_{\mathcal U})=[{\mathcal F}]_{\widetilde{\mathcal U}}$ where
${\mathcal F}$ is a doubly stable filter on $(X,{\mathcal U}).$

$$\begin{xy}
\xymatrix{
(q{\mathcal P}_0(X),q{\mathcal U}_H) \ar[r]^{q\mathcal C} \ar[d]_{qe_0} & (qS_D(X),q{\mathcal U}_D)
 \ar[d]^{qe
 } \\
(q{\mathcal P}_0(\widetilde{X}),q\widetilde{\mathcal U}_H)
\ar[r]_{q\widetilde{\mathcal C}} & (qS_D(\widetilde{X}),q\widetilde{\mathcal U}_D)
}
\end{xy}$$

We next reformulate the K\"unzi-Ryser condition.

\begin{Prp} \label{tres} Let $(X,{\mathcal U})$ be a quasi-uniform
space. Then the Hausdorff quasi-uniformity ${\mathcal U}_H$
is bicomplete on ${\mathcal P}_0(X)$ if and only if each doubly stable filter ${\mathcal F}$
on $(X,{\mathcal U})$ is ${\mathcal U}_D$-equivalent to
some ${\mathcal C}_{C},$ that is
${\mathcal F}_{\mathcal U}={\mathcal D}_{\mathcal U}({\mathcal C}_C)$
for some $C\in {\mathcal P}_0(X).$ (The condition implies that
each doubly stable 2-round filter ${\mathcal F}$ on $(X,{\mathcal U})$
is uniquely determined by some doubly closed set $C\in {\mathcal P}_0(X).$)
\end{Prp}

{\em Proof.}
Suppose that ${\mathcal U}_H$ is bicomplete. Then the
K\"unzi-Ryser condition \cite[Proposition 8]{KunziRyser} is satisfied. It follows that
${\mathcal D}_{\mathcal U}({\mathcal C}_{C({\mathcal F})})$ is coarser than
${\mathcal F}$, where $C({\mathcal F})$ denotes the set of double
cluster points of ${\mathcal F}.$
Clearly by definition of $C({\mathcal F})$
the filter ${\mathcal D}_{\mathcal U}({\mathcal F})$
is coarser than ${\mathcal D}_{\mathcal U}({\mathcal C}_{C({\mathcal F})}).$
Hence ${\mathcal C}_{C({\mathcal F})}$ and ${\mathcal F}$  are
indeed  ${\mathcal U}_D$-equivalent.

In order to prove the converse suppose that for any doubly stable
filter ${\mathcal F}$ on $(X,{\mathcal U})$ there is $C\in
{\mathcal P}_0(X)$ such that ${\mathcal F}$ is ${\mathcal
U}_D$-equivalent to ${\mathcal C}_C.$ Thus ${\mathcal F}_{\mathcal
U}= {\mathcal D}_{\mathcal U}({\mathcal C}_C).$ We immediately
deduce that $C({\mathcal F}):=\mbox{cl}_{\tau({\mathcal U})}C\cap
\mbox{cl}_{\tau({\mathcal U}^{-1})}C$ is equal to the set of
double cluster points of ${\mathcal F}$ in $(X,{\mathcal U}).$
Furthermore ${\mathcal D}_{\mathcal U}({\mathcal C}_C)$ and
${\mathcal D}_{\mathcal U}({\mathcal C}_{C({\mathcal F})})$ are
equal. We conclude that ${\mathcal F}$ is finer than ${\mathcal
D}_{\mathcal U} ({\mathcal C}_{C({\mathcal F})})$, which means
that the K\"unzi-Ryser condition is satisfied.
 $\hfill{\Box}$

\bigskip

The preceding proposition motivates the following results, which deal
with the general case.
They will lead to a characterization of the stability space of a quasi-uniform space.

\begin{Lem} \label{major} Let $(X,{\mathcal U})$ be a quasi-uniform space.

(a) If ${\mathcal F}$ is a doubly stable filter on $(X,{\mathcal U}),$
then $({\mathcal C}_{F})_{F\in ({\mathcal F},\supseteq)}$ is a $({\mathcal U}_D)^s$-Cauchy
net in $(S_D(X),{\mathcal U}_D).$

(b) Let $({\mathcal C}_{C_d})_{d\in E}$ be a $({\mathcal
U}_D)^s$-Cauchy net in the subset ${\mathcal C}{\mathcal P}_0(X)$
of $(S_D(X),{\mathcal U}_D).$ Then there is ${\mathcal F}\in
S_D(X)$ such that $({\mathcal C}_{C_d})_{d\in E}$
$\tau(({{\mathcal U}}_D)^s)$-converges to ${\mathcal F}.$ (For the
net considered in part (a), the constructed point ${\mathcal F}\in
S_D(X)$ is equal to the original filter ${\mathcal F}.$ Therefore
${\mathcal C}{\mathcal P}_0(X)$ is $\tau(({\mathcal
U}_D)^s)$-dense in $S_D(X).$)

(c) The quasi-uniform space $(S_D(X),{\mathcal U}_D)$ is bicomplete.
\end{Lem}

{\em Proof.}
For the convenience of the reader we present a complete proof of
this result, although many techniques are known (see
e.g. \cite{KunziRyser}).

(a) Let ${\mathcal F}$ be a doubly stable filter on
$(X,{\mathcal U})$ and let $U\in {\mathcal U}.$
Then there is $F_U\in {\mathcal F}$ such that $F_U\subseteq U(F)\cap
U^{-1}(F)$
whenever $F\in {\mathcal F}.$
Of course,
$F\subseteq F_U$ implies that
$F\subseteq U(F_U)\cap U^{-1}(F_U).$
Therefore $(F_U,F)\in (U_H)^{-1}\cap U_H$ whenever $F\in {\mathcal F}$ and $F\subseteq F_U.$
By Remark \ref{categ} we have shown that $({\mathcal C}_F)_{F\in ({\mathcal F},\supseteq)}$ is a
$({\mathcal U}_D)^s$-Cauchy net of $(S_D(X),{\mathcal U}_D).$

(b)
Let $({\mathcal C}_{C_d})_{d\in E}$ be a $({\mathcal U}_D)^s$-Cauchy net in the subspace
${\mathcal C}
{\mathcal P}_0(X)$
of $(S_D(X),{{\mathcal U}}_D):$
Therefore for each $U\in {\mathcal U}$ there is $d_U\in E$ such that
for any $d_1,d_2\in E$ satisfying $d_1,d_2\geq d_U$ we have $C_{d_2}\subseteq
U(C_{d_1})\cap U^{-1}(C_{d_1})$ and $C_{d_1}\subseteq U(C_{d_2})\cap U^{-1}(C_{d_2}).$

For each $d\in E$ set $F_d=\bigcup_{d'\in E;d'\geq d}C_{d'}.$
Let ${\mathcal F}$ be the filter generated by the base $\{F_d:d\in E\}$
on $X.$
Let $x\in F_{d_U}$ and $d_2\in E.$ Then $x\in C_{d_1}$ for some $d_1\geq d_U.$
By directedness of $E$ we find $d_3\in E$ such that $d_3\geq d_2,d_U.$ Thus $C_{d_1}\subseteq
U(C_{d_3})\subseteq
U(F_{d_2}).$ Therefore $F_{d_U}\subseteq \bigcap_{F\in {\mathcal F}}U(F).$
Hence ${\mathcal F}$ is ${\mathcal U}$-stable. Similarly it
is shown that ${\mathcal F}$ is ${\mathcal U}^{-1}$-stable. Consequently
${\mathcal F}$ is a doubly stable filter on $(X,{\mathcal U}).$

We further check that the net $({\mathcal C}_{C_d})_{d\in E}$
$\tau(({\mathcal U}_D)^s)$-converges to the constructed point
   ${\mathcal F}$ in
$(S_D(X),{\mathcal U}_D):$
Let $U\in {\mathcal U}.$
Consider any $d\in E$ such that $d\geq d_U.$
As shown above, $F_{d_U}\subseteq \bigcap_{F\in {\mathcal F}}U(F)$
and thus $\bigcap_{F\in {\mathcal F}}U(F)\in {\mathcal C}_{C_d}$
and $({\mathcal F},{\mathcal C}_{C_d})\in U_+.$
We also
have
$F_{d_U}\subseteq U^{-1}(C_d),$ which means
that $U^{-1}(C_d)\in {\mathcal F}$ and $({\mathcal F},{\mathcal
C}_{C_d})\in U_-.$
Similarly,
$\bigcap_{F\in {\mathcal F}}U^{-1}(F)\in {\mathcal C}_{C_d},$
and  $U(C_d)\in {\mathcal F}.$ Hence $({\mathcal C}_{C_d},{\mathcal F})\in
U_D$ and the claim is verified.
We conclude that the net
$({\mathcal C}_{C_d})_{d\in E}$ $\tau(({\mathcal U}_D)^s)$-converges
to ${\mathcal F}$
in $(S_D(X),{\mathcal U}_D).$

The final assertion about ${\mathcal F}$ is obvious by the construction of the filter ${\mathcal F}.$
Hence if we start with a doubly stable filter ${\mathcal F}$
on $(X,{\mathcal U}),$ the presented argument shows
that the net
$({\mathcal C}_F)_{F\in ({\mathcal F},\supseteq)}$
$\tau(({\mathcal U}_D)^s)$-converges to the point ${\mathcal F}$
in $(S_D(X),{\mathcal U}_D).$
In particular this proof establishes
that ${\mathcal C}{\mathcal P}_0(X)$ is $\tau(({\mathcal U}_D)^s)$-dense in
$(S_D(X),{{\mathcal U}}_D).$

(c) We shall give a direct proof that $(S_D(X),{{\mathcal U}}_D)$
is bicomplete (but compare \cite[Proposition 3.32] {FletcherLindgren}).
Let ${\mathcal F}$ be a $({\mathcal U}_D)^s$-Cauchy filter on
$(S_D(X),{\mathcal U}_D).$
Thus for each $U\in {\mathcal U},$ there is
$F_U\in {\mathcal F}$ such that
   $F_U\times F_U\subseteq {U}_D.$
For each $U\in {\mathcal U}$ we find
${\mathcal C}_{C_U}\in {\mathcal C}{\mathcal P}_0(X)\cap
({U}_D)^{-1}(F_U)\cap {U}_D(F_U)$ by
$\tau(({\mathcal U}_D)^s)$-density of ${\mathcal C}{\mathcal P}_0(X)$
in $S_D(X).$

Then $({\mathcal C}_{C_U})_{U\in ({\mathcal U},\supseteq)}$
   is a $({\mathcal U}_D)^s$-Cauchy net on the subspace ${\mathcal C}{\mathcal P}_0(X)$ of
   $(S_D(X),{{\mathcal U}}_D):$
   Let $U\in {\mathcal U}.$ Choose $V\in {\mathcal U}$ such that
${V}^5\subseteq {U}.$ Thus $({V}_D)^5\subseteq {U}_D$ (see proof
of Proposition \ref{stabcon}).
Let $P_1,P_2\in {\mathcal U}$ be such that $P_1,P_2\subseteq V.$
Then ${\mathcal C}_{C_{P_1}}\in (({P_1})_D)^{-1}({\mathcal F}_{P_1})$ where
${\mathcal F}_{P_1}\in F_{P_1}$ and $F_{P_1}\times F_{P_1}\subseteq
{(P_1)}_D.$
Furthermore ${\mathcal C}_{C_{P_2}}\in ({P_2})_D({\mathcal F}_{P_2})$ where
${\mathcal F}_{P_2}\in F_{P_2}$ and $F_{P_2}\times F_{P_2}\subseteq (P_2)_D.$
We find ${\mathcal A}\in F_{P_1}\cap F_V$ and ${\mathcal B}\in F_{P_2}\cap F_V.$
Consequently $({\mathcal C}_{C_{P_1}},{\mathcal F}_{P_1})\in ({P}_1)_D,$
$({\mathcal F}_{P_1},{\mathcal A})\in F_{P_1}\times F_{P_1}\subseteq
{V}_D,$ $({\mathcal A},{\mathcal B})
\in F_V\times F_V\subseteq {V}_D,
({\mathcal B},{\mathcal F}_{P_2})\in F_{P_2}\times F_{P_2}\subseteq V_D$
and $({\mathcal F}_{P_2},{\mathcal C}_{C_{P_2}})\in ({{P}_2})_D.$
Thus $({\mathcal C}_{C_{P_1}},{\mathcal C}_{C_{P_2}})\in
({V}_D)^5\subseteq {U}_D.$

We have proved  that
$({\mathcal C}_{C_U})_{U \in ({\mathcal U},\supseteq)}$ is a
$({\mathcal{U}}_D)^s$-Cauchy net in ${\mathcal C}{\mathcal P}_0(X).$
   Thus according to the argument above there exists a doubly stable
filter ${\mathcal H}$ on
   $(X,{\mathcal U})$ to which it
   $\tau(({\mathcal U}_D)^s)$-converges
   in $(S_D(X),{\mathcal U}_D).$
   Clearly then ${\mathcal F}$ also
   $\tau(({\mathcal U}_D)^s)$-converges to ${\mathcal H}:$
   Indeed given $U\in {\mathcal U}$ there is $P\in {\mathcal U}$ such
that $P\subseteq U$ and
   $(\mathcal C_{C_P},{\mathcal H})\in {U}_D.$
   By definition of ${\mathcal C}_{C_P}$ we find ${\mathcal F}_P\in
F_P$ such that $({\mathcal F}_P,{\mathcal C}_{C_P})\in {P}_D.$
   Recall that $F_{P}\times F_P\subseteq {P}_D.$ Let ${\mathcal B}_P\in F_P.$
   Then we have
$({\mathcal B}_P,{\mathcal F}_P)\in {P}_D.$
Therefore $({\mathcal B}_P,{\mathcal F}_P)\in U_D,$
$({\mathcal F}_P,{\mathcal C}_{C_P})\in U_D$ and
$({\mathcal C}_{C_P},{\mathcal H})\in U_D.$
Thus $F_P\subseteq ({U}_D)^{-3}({\mathcal H}).$

Convergence of ${\mathcal F}$ to ${\mathcal H}$ in the conjugate
topology $\tau({\mathcal U}_D)$ is established analogously.
We finally conclude that the filter
${\mathcal F}$ $\tau(({\mathcal U}_D)^s)$-converges to ${\mathcal H}$ in
$(S_D(X),{\mathcal U}_D).$
$\hfill{\Box}$

\begin{Thm} \label{main} Let $(X,{\mathcal U})$ be a quasi-uniform space. Then the
$T_0$-quotient $(qS_D(X),q{\mathcal U}_D)$ of $(S_D(X),{\mathcal
U}_D)$ is the bicompletion of the subspace $q{\mathcal C}{\mathcal
P}_0(X):=\{{\mathcal D}_{\mathcal U}({\mathcal C}_C):C\in
{\mathcal P}_0(X), C$ is doubly closed in $(X,{\mathcal U})\}.$
(We note that according to Lemma \ref{rep} $qS_D(X)$ can be
identified with the set of all doubly stable 2-round filters on
$(X,{\mathcal U}).$)
\end{Thm}

{\em Proof.} The statement is a consequence of the preceding
result, since it is well known (and easy to see) that a quasi-uniform space is
bicomplete if and only if its $T_0$-quotient is bicomplete.
From Lemma \ref{major} it follows that $(qS_D(X),q{\mathcal U}_D)$ is bicomplete. Furthermore
$q{\mathcal C}{\mathcal P}_0(X)$ is $\tau((q{\mathcal U}_D)^s)$-dense in $qS_D(X).$
For this we note similarly as above that for any ${\mathcal F}\in qS_D(X)$
we have that the net $$({\mathcal D}_{\mathcal U}({\mathcal C}_{F'}))_{(F'\in {\mathcal F},
F' {\small\mbox{ is doubly closed in }} (X,{\mathcal U}),\,\supseteq)}$$ $\tau((q{\mathcal U}_D)^s)$-converges to
${\mathcal F}.$ Hence the assertion is proved.
$\hfill{\Box}$

\begin{Prp} \label{bc} Let $(X,{\mathcal U})$ be a quasi-uniform $T_0$-space
and $(\widetilde{X},\widetilde{\mathcal U})$ its bicompletion.
Then the $T_0$-quotients $(qS_D(X),q{\mathcal U}_D)$
of $(S_D(X),{\mathcal U}_D)$ and
$(qS_D(\widetilde{X}),q\widetilde{{\mathcal U}}_D)$ of
$(S_{D}(\widetilde{X}),\widetilde{{\mathcal U}}_D)$
are isomorphic as quasi-uniform spaces under the quasi-uniform isomorphism $qe$
(see the second diagram above).

Hence by Theorem \ref{main} $(qS_D(X),q{\mathcal U}_D)$ can also be understood
as the bicompletion of the image of
$(q{\mathcal P}_0(\widetilde{X}),q\widetilde{{\mathcal U}}_H)$
under the quasi-uniform embedding $q\widetilde{\mathcal C}$
into $(qS_D(\widetilde{X}),q\widetilde{{\mathcal U}}_D).$
\end{Prp}

{\em Proof.} Recall that $qe:qS_D(X)\rightarrow
qS_D(\widetilde{X})$ is defined by $qe(({\mathcal F}')_{\mathcal
U})=[{\mathcal F}']_{\widetilde{\mathcal U}},$ where ${\mathcal
F}'$ is a doubly stable filter on $(X,{\mathcal U}).$ Then
$[{\mathcal F}']_{\widetilde{\mathcal U}}={[\mathcal
G']}_{\widetilde{\mathcal U}}$ (with ${\mathcal F}',$ ${\mathcal
G}'\in S_D(X))$ clearly implies that ${({\mathcal F}')}_{\mathcal
U}={({\mathcal G}')}_{\mathcal U}.$ Thus $qe$ is injective.

Suppose that ${\mathcal F}\in S_D(\widetilde{X}).$
Then $qe(({\mathcal F}_{\widetilde{\mathcal U}}\vert X)_{\mathcal U})=
[{\mathcal F}_{\widetilde{\mathcal U}}\vert X]_{\widetilde{\mathcal U}}=
{\mathcal F}_{\widetilde{\mathcal U}}$ by Corollary \ref{resbic}.
Thus $qe$ is surjective. For later use observe that $(qe)^{-1}({\mathcal F}_{\widetilde{\mathcal U}})
=({\mathcal F}_{\widetilde{\mathcal U}}\vert X)_{\mathcal U}.$

We next give a proof from first principles that $qe$ is quasi-uniformly continuous.
Let $U\in {\mathcal U}.$ Choose $V\in {\mathcal U}$ such that $V^2\subseteq U.$
Without loss of generality
\footnote{Here we recall that each quasi-uniformity ${\mathcal U}$ has a base consisting of
$\tau({\mathcal U}^{-1})\times \tau({\mathcal U})$-open entourages $V$
\cite[Corollary 1.17]{FletcherLindgren}.
 Such entourages $V$ obviously satisfy the equality
$\widetilde{V}\cap (X\times X) =V.$}
we can assume that $V$ is
$\tau({\mathcal U}^{-1})\times \tau({\mathcal U})$-open.

Furthermore let ${\mathcal F'},{\mathcal G'}\in S_D(X)$ be such that
$(({\mathcal F'})_{\mathcal U},({\mathcal G'})_{\mathcal U})\in V_D.$
Therefore there is $G\in ({\mathcal G'})_{\mathcal U}$
such that $G\subseteq V(F)$ whenever $F\in ({\mathcal F'})_{\mathcal U}.$
Thus $\widetilde{V}(G)\subseteq \widetilde{V}^2(F)\subseteq \widetilde{U}(F)$
whenever $F\in ({\mathcal F'})_{\mathcal U},$
and therefore $\widetilde{V}(G)\subseteq \widetilde{U}(E)$ whenever
$E\in [{\mathcal F'}]_{\widetilde{\mathcal U}}.$
Consequently $\bigcap_{E\in [{\mathcal F'}]_{\widetilde{\mathcal U}}}
\widetilde{U}(E)\in [{\mathcal G'}]_{\widetilde{\mathcal U}}.$
Similarly $\bigcap_{H\in [{\mathcal G'}]_{\widetilde{\mathcal U}}}
\widetilde{U}^{-1}(H)\in [{\mathcal F'}]_{\widetilde{\mathcal U}}.$
Thus $([{\mathcal F'}]_{\widetilde{\mathcal U}},[{\mathcal G'}]_{\widetilde{\mathcal U}})\in
\widetilde{U}_D.$ We conclude that $qe$ is quasi-uniformly continuous.

 It remains
to show that $(qe)^{-1}$ is quasi-uniformly continuous, too.
Let $U\in {\mathcal U}.$ There are $W,V\in {\mathcal U}$ such that $W^2\subseteq V$
and $V^2\subseteq U.$ Again we assume that
$V$ is
$\tau({\mathcal U}^{-1})\times \tau({\mathcal U})$-open.

Let ${\mathcal F},{\mathcal G}\in S_D(\widetilde{X})$
such that $({\mathcal F}_{\widetilde{\mathcal U}},{\mathcal G}_{\widetilde{\mathcal U}})\in
\widetilde{W}_D.$
We want to show that
$(({\mathcal F}_{\widetilde{\mathcal U}}\vert X)_{\mathcal U},
({\mathcal G}_{\widetilde{\mathcal U}}\vert X)_{\mathcal U})\in
U_D.$
There is $G\in {\mathcal G}_{\widetilde{\mathcal U}}$ such that
$G\subseteq \widetilde{W}(F)$ whenever $F\in {\mathcal F}_{\widetilde{\mathcal U}}\cap
\tau(\widetilde{\mathcal U}^s).$
For such $F$ we also have
 $\widetilde{W}(F)\subseteq \widetilde{W}(\mbox{cl}_{\tau(\widetilde{\mathcal U}^s)}(F\cap X))
\subseteq \widetilde{W}^2(F\cap X)\subseteq \widetilde{V}(F\cap
X).$ Consequently $G\cap X\subseteq V(F\cap X)$ and hence $V(G\cap
X)\subseteq V^2(F\cap X)\subseteq U(Z(F\cap X))\cap U(Z^{-1}(F\cap
X))$ whenever $Z\in {\mathcal U}$ and $F\in {\mathcal
F}_{\widetilde{{\mathcal U}}}.$ It follows that $\bigcap_{E\in
({\mathcal F}_{\widetilde{\mathcal U}} \vert X)_{\mathcal
U}}U(E)\in ({\mathcal G}_{\widetilde{\mathcal U}}\vert
X)_{\mathcal U}.$ Similarly $\bigcap_{G\in ({\mathcal
G}_{\widetilde{\mathcal U}}\vert X)_{\mathcal U}}U^{-1}(G) \in
({\mathcal F}_{\widetilde{\mathcal U}}\vert X)_{\mathcal U}.$ Thus
$(({\mathcal F}_{\widetilde{\mathcal U}}\vert X)_{\mathcal U},
({\mathcal G}_{\widetilde{\mathcal U}}\vert X)_{\mathcal U})\in
U_D$ and we are finished. $\hfill{\Box}$

\bigskip

We next want to address the problem of characterizing those
quasi-uniform $T_0$-spaces $(X,{\mathcal U})$
such that $({\mathcal
P}_0(\widetilde{X}),
\widetilde{{\mathcal U}}_H)$ is bicomplete, where $(\widetilde{X},\widetilde{\mathcal U})$
is the bicompletion of $(X,{\mathcal U}).$

\begin{Rem} \label{st} It is straightforward to check that on a quasi-uniform space
$(X,{\mathcal U})$ the intersection of each nonempty family
of (doubly) stable filters is (doubly) stable. In particular the intersection
of any nonempty family of ${\mathcal U}^s$-Cauchy filters on a quasi-uniform
space $(X,{\mathcal U})$ is doubly stable, since such filters are ${\mathcal U}^s$-stable.
These observations motivated the following investigations.
\end{Rem}

\begin{Lem} \label{aux} Let $(X,{\mathcal U})$ be a quasi-uniform $T_0$-space and
$(\widetilde{X},\widetilde{\mathcal U})$ its bicompletion.
For each $C\in {\mathcal P}_0(\widetilde{X})$ we have that
$(\widetilde{C}_C)_{\widetilde{\mathcal U}}\vert X$ and
$\bigcap_{x\in C}(\widetilde{\mathcal U}^s(x)\vert X)$
are ${\mathcal U}_D$-equivalent filters on $X.$

\end{Lem}

{\em Proof.} By Remark \ref{st} and Lemma \ref{dense} it is obvious that the two
filters under consideration are doubly stable on $(X,{\mathcal U}).$
We have that $\bigcap_{x\in C}\widetilde{\mathcal U}^s(x)
\subseteq {\widetilde{\mathcal C}}_C,$ since $C\subseteq A$
whenever $A\in \bigcap_{x\in C}\widetilde{\mathcal U}^s(x)$.

Therefore
$$(\bigcap_{x\in C}\widetilde{\mathcal U}^s(x))_{\widetilde{\mathcal U}}
\subseteq
({\widetilde{\mathcal C}}_C)_{\widetilde{\mathcal U}}.$$

Let $U\in {\mathcal U}.$
For each $x\in C\subseteq \widetilde{X}$
consider $\widehat{U}_x
\in {\widetilde{\mathcal U}}$ such that
${\widehat{U}_x}$ is $\tau(\widetilde{{\mathcal U}}^{-1})\times
\tau(\widetilde{{\mathcal U}})$-open (see \cite[Corollary 1.17]{FletcherLindgren}).
Note that then $\widehat{U}^s_x(x)$ is $\tau(\widetilde{{\mathcal U}}^s)$-open.

It follows that

$$\widetilde{U}(\bigcup_{x\in C}\widehat{U}^s_x(x))\subseteq \widetilde{U}
(\mbox{cl}_{\tau(\widetilde{\mathcal U}^s)}(\bigcup_{x\in C}\widehat{U}^s_x(x)))$$
$$\subseteq \widetilde{U}(\mbox{cl}_{\tau(\widetilde{\mathcal U}^s)} ((\bigcup_{x\in C}
\widehat{U}^s_x(x))\cap X))\subseteq \widetilde{U}
(\mbox{cl}_{\tau(\widetilde{\mathcal U}^s)}\bigcup_{x\in C}(
\widehat{U}^s_x(x)\cap X))\subseteq \widetilde{U}^2(\bigcup_{x\in C}(
\widehat{U}^s_x(x)\cap X)).$$

The conjugate inequality for $\widetilde{U}^{-1}$ is established similarly.

We deduce by Lemma \ref{drei} that

$$[\bigcap_{x\in C}(\widetilde{\mathcal U}^s(x)\vert X)]_{\widetilde{\mathcal U}}
\subseteq (\bigcap_{x\in C}\widetilde{\mathcal
U}^s(x))_{\widetilde{\mathcal U}}\subseteq ({\widetilde{\mathcal
C}}_C)_{\widetilde{\mathcal U}}.$$

Furthermore $({\widetilde{\mathcal C}}_{C})_{\widetilde{\mathcal U}}\subseteq
[\bigcap_{x\in C}
(\widetilde{\mathcal U}^s(x)\vert X)],$ because $\bigcup_{x\in C} (\widetilde{U}^s(x)\cap X)\subseteq
\widetilde{U}(C)\cap \widetilde{U}^{-1}(C)$ whenever $U\in {\mathcal U}.$
It follows that $(\widetilde{{\mathcal C}}_{C})_{\widetilde{\mathcal U}}\subseteq [\bigcap_{x\in C}
(\widetilde{\mathcal U}^s(x)\vert X)]_{\widetilde{\mathcal U}}.$
Altogether therefore $[\bigcap_{x\in C}(\widetilde{\mathcal U}^s(x)\vert X)]_{\widetilde{\mathcal U}}
=({\widetilde{\mathcal C}_C)}_{\widetilde{\mathcal U}}=
[({\widetilde{\mathcal C}_C)}_{\widetilde{\mathcal U}}\vert X]_{\widetilde{\mathcal U}}$
by Corollary \ref{resbic}.

We conclude that
$$(\bigcap_{x\in C}(\widetilde{\mathcal U}^s(x)\vert X))_{{\mathcal U}}
=(({\widetilde{\mathcal C}_C})_{\widetilde{\mathcal U}}\vert X)_{{\mathcal U}}.$$
$\hfill{\Box}$

\begin{Prp} \label{bchar} Let $(X,{\mathcal U})$ be a quasi-uniform
$T_0$-space and
 $(\widetilde{X},\widetilde{\mathcal U})$
 its bicompletion.
Then $({\mathcal P}_0(\widetilde{X}), \widetilde{\mathcal U}_H)$
is bicomplete if and only if each doubly stable filter on
$(X,{\mathcal U})$ is ${\mathcal U}_D$-equivalent to the
intersection of a nonempty family of ${\mathcal U}^s$-Cauchy
filters on $(X,{\mathcal U}).$
\end{Prp}

{\em Proof.} Suppose that $({\mathcal P}_0(\widetilde{X}),
\widetilde{\mathcal U}_H)$ is bicomplete and let ${\mathcal F}$ be
a doubly stable filter on $(X,{\mathcal U}).$ Then by the proof of
Proposition \ref{tres} $[{\mathcal F}]_{\widetilde{\mathcal U}}
=({\widetilde{\mathcal C}}_C)_{\widetilde{\mathcal U}}$ where $C$
is the nonempty set of double cluster points of $[{\mathcal F}]$
in $(\widetilde{X},\widetilde{\mathcal U}).$

We show that ${\mathcal F}$ and $\bigcap_{x\in
C}(\widetilde{\mathcal U}^s(x)\vert X)$ are ${\mathcal
U}_D$-equivalent: By Lemma \ref{aux} and Corollary \ref{resbic}
${\mathcal F}_{\mathcal U}= ([{\mathcal F}]_{\widetilde{\mathcal
U}}\vert X)_{\mathcal U}= ((\widetilde{\mathcal
C}_C)_{\widetilde{\mathcal U}}\vert X)_{\mathcal U}
=(\bigcap_{x\in C}(\widetilde{\mathcal U}^s(x)\vert X))_{\mathcal
U}.$ Thus ${\mathcal F}$ is ${\mathcal U}_D$-equivalent to
$\bigcap_{x\in C}(\widetilde{\mathcal U}^s(x)\vert X)$ on $X$
where each $\widetilde{\mathcal U}^s(x)\vert X$ is a ${\mathcal
U}^s$-Cauchy filter on $X.$

For the converse suppose that ${\mathcal F}$ is a doubly stable
filter on $(\widetilde{X},\widetilde{\mathcal U}).$ Set ${\mathcal
F}'={\mathcal F}_{\widetilde{\mathcal U}}\vert X.$ Then by Lemma
\ref{dense} ${\mathcal F}'$ is a doubly stable filter on
$(X,{\mathcal U})$ that by our assumption is ${\mathcal
U}_D$-equivalent to $\bigcap_{i\in I}{\mathcal F}'_i$ where $I\not
= \emptyset$ and each ${\mathcal F}'_i$ is a ${\mathcal
U}^s$-Cauchy filter on $(X,{\mathcal U}).$ Let $C=\{x\in
\widetilde{X}:$ there is $i\in I$ such that $[{\mathcal F}_i']$
$\tau(\widetilde{{\mathcal U}}^s)$-converges to $x\}.$\footnote{It
is well known (and easy to see) that in a
 quasi-uniform $T_0$-space $(X,{\mathcal V})$
$\tau({\mathcal V}^s)$-limits of filters are unique.}

We wish to show that ${\mathcal F}_{\widetilde{\mathcal U}}
=(\widetilde{{\mathcal C}}_C) _{\widetilde{\mathcal U}}.$ Let
$U\in {\mathcal U}$ and $F'\in {\mathcal F}'.$ By assumption for
each $i\in I$ there is $F'_i\in {\mathcal F}_i'$ such that
$F'_i\subseteq \widetilde{U}(F')\cap \widetilde{U}^{-1}(F').$
Suppose that for each $i\in I,$ $x_i$ denotes the
$\tau(\widetilde{\mathcal U}^s)$-limit of $[{\mathcal F}_i']$ on
$\widetilde{X}.$ Then for each $i\in I,$ $x_i\in
\widetilde{U}(F'_i)\cap {\widetilde{U}}^{-1}(F'_i)$ and
consequently $x_i\in \widetilde{U}^{2}(F')\cap
\widetilde{U}^{-2}(F').$ Thus $[{\mathcal
F'}]_{\widetilde{\mathcal U}} \subseteq \widetilde{{\mathcal
C}}_C$ and $[{\mathcal F}']_{\widetilde{\mathcal U}}
\subseteq(\widetilde{{\mathcal C}}_C)_{\widetilde{\mathcal U}}.$

On the other hand given $U\in {\mathcal U},$ we have
$\bigcup_{x\in C}(\widetilde{U}^s(x)\cap X)\subseteq
\widetilde{U}(C)\cap \widetilde{U}^{-1}(C).$ Note that
$\bigcup_{x\in C}(\widetilde{U}^s(x)\cap X)\in \bigcap_{i\in
I}{\mathcal F}_i'$ since for each $i\in I,$ $[F_i']$
$\tau(\widetilde{\mathcal U}^s)$-converges to some $x\in C.$

Thus $({\widetilde{\mathcal C}}_C)_{\widetilde{\mathcal U}}=
{\mathcal D}_{\widetilde{\mathcal U}}({\widetilde{\mathcal
C}}_C)\subseteq [\bigcap_{i\in I}{\mathcal F}_i']$ and therefore
$(\widetilde{{\mathcal C}}_C) _{\widetilde{\mathcal U}}\subseteq
[\bigcap_{i\in I}{\mathcal F}_i'] _{\widetilde{\mathcal U}}
=[{\mathcal F}']_ {\widetilde{\mathcal U}}.$ We conclude that
$({\widetilde{\mathcal C}}_C) _{\widetilde{\mathcal
U}}=[{{\mathcal F}'}] _{\widetilde{\mathcal U}}.$

Since by Corollary \ref{resbic} $[{\mathcal F}']_{\widetilde{\mathcal U}}=
{\mathcal F}_{\widetilde{\mathcal U}},$
it follows that $({\mathcal P}_0(\widetilde{X}),\widetilde{\mathcal U}_H)$
is bicomplete by Proposition \ref{tres}.
$\hfill{\Box}$

\bigskip
Observe by the preceding result that in a complete uniform
$T_0$-space $(X,{\mathcal U})$ which is not supercomplete there
must exist a stable filter that is not ${\mathcal U}_D$-equivalent
to the intersection of a nonempty family of ${\mathcal U}$-Cauchy
filters on $X.$

\begin{Rem}
We note that the techniques to establish Propositions \ref{tres} and \ref{bchar} can be
combined to yield the following result:
 For a quasi-uniform $(X,\mathcal{U})$ the quasi-uniform space
 $(\mathcal{P}_0(X), \mathcal{U}_H)$ is bicomplete if and only if each
 doubly stable 2-round filter on $(X,\mathcal{U})$ is the
 2-envelope of the intersection of a nonempty family of
 $\mathcal{U}^s$-convergent filters on $X.$
 $\hfill{\Box}$
\end{Rem}

We remarked above that the Burdick-Isbell criterion for
completeness of the Hausdorff uniformity is easier to understand than the K\"unzi-Ryser condition that characterizes bicompleteness
of the Hausdorff quasi-uniformity. We next show that analogously Proposition
\ref{bchar} can be simplified in the case that we are only interested in
uniform spaces.

\begin{Prp} \label{unifor} Let $(X,{\mathcal U})$ be a uniform $T_0$-space
and $(\widetilde{X},\widetilde{\mathcal U})$ its completion. Then
$({\mathcal P}_0(\widetilde{X}),{\widetilde{\mathcal U}}_H)$ is complete
if and only if each stable filter on $(X,{\mathcal U})$
is contained in a ${\mathcal U}$-Cauchy filter.
\end{Prp}

{\em Proof.} Suppose that $({\mathcal P}_0(\widetilde{X}),\widetilde{\mathcal U}_H)$
is complete and let ${\mathcal F}$ be a stable filter on $(X,{\mathcal U}).$
Then by Proposition \ref{bchar} ${\mathcal F}$ is
 ${\mathcal U}_D$-equivalent to the intersection of a
nonempty family $({\mathcal F}_i)_{i\in I}$
of ${\mathcal U}$-Cauchy filters on $X.$ So
${\mathcal F}_{\mathcal U}=(\bigcap_{i\in I}{\mathcal F}_i)_{\mathcal U}.$
For any ${\mathcal U}$-Cauchy filter
${\mathcal F}'$ of this family we have that
$U(F)\cap F'\not = \emptyset$ whenever $U\in {\mathcal U},$
$F \in {\mathcal F}$
and $F'\in {\mathcal F}'.$ We conclude that the filterbase
$\{F\cap U(F'):F\in {\mathcal F}, F' \in{\mathcal F}',U\in {\mathcal U}\}$
generates
a ${\mathcal U}$-Cauchy filter on $X$ that is finer than ${\mathcal F}.$
Hence the stated criterion is satisfied.

For the converse suppose that each stable filter ${\mathcal F}$
on the uniform $T_0$-space  $(X,{\mathcal U})$ is coarser than a ${\mathcal U}$-Cauchy
filter on $X.$
Consider now an arbitrary stable filter ${\mathcal F}$ on $(X,{\mathcal U}).$
Let ${\mathcal H}=\{{\mathcal G}:{\mathcal G}$ is a ${\mathcal U}$-Cauchy filter finer
than ${\mathcal F}$ on $(X,{\mathcal U})\}.$
We want to show that $\bigcap_{{\mathcal G}\in {\mathcal H}}{\mathcal G}$
and ${\mathcal F}$ are ${\mathcal U}_D$-equivalent.
Certainly ${\mathcal F}\subseteq \bigcap_{{\mathcal G}\in
{\mathcal H}}{\mathcal G}$ by definition of ${\mathcal H}$
and thus
${\mathcal F}_{\mathcal U}\subseteq (\bigcap_{{\mathcal G}\in
{\mathcal H}}{\mathcal G})_{\mathcal U}.$

In order to reach a contradiction suppose that $(\bigcap_{{\mathcal G}\in
{\mathcal H}}{\mathcal G})_{\mathcal U}\not \subseteq {\mathcal F}.$
Thus
there is $U_0\in {\mathcal U}$ and
 for each ${\mathcal G}\in {\mathcal H}$ there is
$G_{\mathcal G}\in {\mathcal G}$
such that $F\setminus
U^2_{0}(\bigcup_{{\mathcal G}\in {\mathcal H}}G_{\mathcal G})
\not =\emptyset$ whenever $F\in {\mathcal F}.$
For each $U\in {\mathcal U}$ and $E\in {\mathcal F}$
set $H_{UE}=\{a\in X:$ there is $V\in {\mathcal U}$ such
that $V^2\subseteq U, V^{-2}(a)\cap U_0(\bigcup_{{\mathcal G}\in
{\mathcal H}} G_{\mathcal G})$ is empty and $a\in
\bigcap_{F\in {\mathcal F}}V(F)\cap E\}.$ According to the proof
of \cite[Lemma 6]{KunziRyser} $\{H_{UE}:U\in {\mathcal U},
E\in {\mathcal F}\}$ is a base for a stable filter on $(X,{\mathcal U}).$
Thus it is contained in a ${\mathcal U}$-Cauchy filter ${\mathcal K'}$
on $X$ by
our assumption. Since ${\mathcal K}'\in {\mathcal H},$
 we see that
$G_{\mathcal K'}\in {\mathcal K'},$ as well as $X\setminus G_{\mathcal K'}
\in {\mathcal K'}$ by the definition of the sets $H_{UE}$
---a contradiction. We conclude that
 $(\bigcap_{{\mathcal G}\in
{\mathcal H}}{\mathcal G})_{\mathcal U} \subseteq {\mathcal F}_{\mathcal U}.$
Therefore
 ${\mathcal F}$ and
$\bigcap_{{\mathcal G}\in {\mathcal H}}{\mathcal G}$ are
${\mathcal U}_D$-equivalent. Hence by Proposition \ref{bchar} $({\mathcal P}_0(\widetilde{X}),
\widetilde{{\mathcal U}}_H)$ is complete.
$\hfill{\Box}$

\section{An application of the stability construction}

It is interesting to study the stability quasi-uniformity,
either using direct proofs or known facts
about the Hausdorff quasi-uniformity
and the bicompletion. For either method we present an illustrating example.

\begin{Prp} Let $(X,{\mathcal U})$ be a quasi-uniform space. Then
$(S_D(X),{\mathcal U}_D)$ is precompact if and only if $(X,{\mathcal U})$
is precompact.
\end{Prp}

{\em Proof.} The reader may want to compare the following
argument with the proof of \cite[Proposition 1]{KunziRyser}, where the
analogous result for the Hausdorff quasi-uniformity ${\mathcal U}_H$ is established.

 Let $(X,{\mathcal U})$ be precompact and let $V\in
{\mathcal U}_D.$
There are $W,U\in {\mathcal U}$ such that $W^2\subseteq U$ and
$U_D\subseteq V.$ Since ${\mathcal U}$ is precompact, there exists a finite
set $F\subseteq X$ such that $\bigcup_{f\in F}W(f)=X.$ Set
${\mathcal M}={\mathcal P}_0(F).$
We want to show that $S_D(X)=\bigcup_{E\in {\mathcal M}}U_D({\mathcal C}_E):$
Consider an arbitrary ${\mathcal F}\in S_D(X).$
Set $F_{\mathcal F}=\{f\in F:W_{\mathcal F}\cap W(f)\not=\emptyset\}.$
Thus $F_{\mathcal F}\subseteq W^{-1}(W_{\mathcal F})\subseteq
\bigcap_{F\in {\mathcal F}}W^{-2}(F).$ It follows that $\bigcap_{F\in {\mathcal
F}}U^{-1}(F)\in {\mathcal C}_{F_{\mathcal F}}$
and ${\mathcal F}\in U_-({\mathcal C}_{F_{\mathcal F}}).$

Furthermore ${\mathcal F}\in U_+({\mathcal C}_{F_{\mathcal F}}),$
because
$W_{\mathcal F}\subseteq \bigcup\{W(f):W_{\mathcal F}\cap W(f)\not
=\emptyset\}=W(F_{\mathcal F}),$ and therefore
$W_{\mathcal F}\subseteq U(F_{\mathcal F})$ and thus
$U(F_{\mathcal F})\in {\mathcal F}.$ Therefore ${\mathcal F}\in U_+({\mathcal C}_{F_{\mathcal F}}).$
We conclude that $(S_D(X),{\mathcal U}_D)$ is precompact.

On the other hand, suppose that $(S_D(X),{\mathcal U}_D)$ and thus
$(S_D(X),{\mathcal U}_-)$ is precompact. Let $U, V\in {\mathcal
U}$ be such that $V^2\subseteq U.$ By our assumption there is a
finite subcollection ${\mathcal H}$ of $S_D(X)$ such that for each
${\mathcal F}\in S_D(X)$ there is ${\mathcal A}\in {\mathcal H}$
with $\bigcap_{F\in {\mathcal F}}V^{-1}(F)\in {\mathcal A}.$ For
each ${\mathcal A}\in {\mathcal H}$ choose some $x_{\mathcal A}\in
V_{\mathcal A}.$ Then $B=X\setminus \bigcup_{{\mathcal A}\in
{\mathcal H}}U(x_{\mathcal A})$ is necessarily empty. Otherwise
${\mathcal C}_B\in V_-({\mathcal A})$ for some ${\mathcal A}\in
{\mathcal H}.$ Note that ${\mathcal A}\in V_-({\mathcal
C}_{V_{\mathcal A}}),$ since $V_{\mathcal A}\subseteq
\bigcap_{A\in {\mathcal A}}V^{-1}(A).$ Consequently ${\mathcal
C}_B\in (V_-)^2({\mathcal C}_{V_{\mathcal A}})\subseteq
U_-({\mathcal C}_{V_{\mathcal A}}).$ Thus $V_{\mathcal A}
\subseteq U^{-1}(B)$. But then $U(x_{\mathcal A})\cap B\not
=\emptyset.$ Therefore we have reached a contradiction and
conclude that $(X,{\mathcal U})$ is precompact. $\hfill{\Box}$

\begin{Prp} A quasi-uniform space $(X,{\mathcal U})$ is totally
bounded if and only if $(S_D(X),{\mathcal U}_D)$
is totally bounded.
\end{Prp}

{\em Proof.}
Note that total boundedness is preserved
under subspaces \cite[p. 12]{FletcherLindgren}.
Since by Remark \ref{categ}
 $x\mapsto {\mathcal C}_{\{x\}}$ where $x\in X$ yields an embedding of $(X,{\mathcal U})$
into $(S_D(X),{\mathcal U}_D)$, the space $(X,{\mathcal U})$ is totally bounded if
$(S_D(X),{\mathcal U}_D)$ is totally bounded.

For the converse observe that total boundedness is preserved under the Hausdorff
hyperspace construction \cite[Corollary 2]{KunziRyser}
 as well as under the bicompletion \cite[Proposition 3.36]{FletcherLindgren}.
Furthermore a quasi-uniform
space is totally bounded if and only if its $T_0$-quotient is totally bounded.
We conclude by Theorem \ref{main} that $(S_D(X),{\mathcal U}_D)$ is totally bounded if
$(X,{\mathcal U})$ is totally bounded.
$\hfill{\Box}$

{\footnotesize

\noindent
Hans-Peter A. K\"unzi\\
Department of Mathematics and Applied Mathematics\\
University of Cape Town\\
Rondebosch 7701\\
South Africa

\noindent
E-mail: hans-peter.kunzi@uct.ac.za

\bigskip \noindent
S. Romaguera\\
Instituto Universitario de Matem\'atica Pura y Aplicada,\\
Universidad Polit\'ecnica de Valencia\\
46071 Valencia\\
Spain

\noindent
E-mail: sromague@mat.upv.es

\bigskip \noindent
M.A. S\'anchez Granero \\
Area of Geometry and Topology\\
Faculty of Science \\
Universidad de Almer\'{\i}a\\
Spain

\noindent
E-mail: misanche@ual.es}


\begin{thebibliography}{XYZ}


\bibitem{ArticoMarconiPelant} G. Artico, U. Marconi and J. Pelant,
{\em On supercomplete $\omega_{\mu}$-metric spaces},
Bull. Polish Acad. Sci. Math. {\bf 44} (1996), 299--310.

\bibitem{Berthiaume} G. Berthiaume, {\em On
quasi-uniformities in hyperspaces}, Proc. Amer. Math.
Soc. {\bf 66} (1977), 335--343.

\bibitem{BuhagiarPasynkov} D. Buhagiar and B.A. Pasynkov,
{\em Supercomplete topological spaces}, Acta Math. Hungar.
{\bf 115(4)} (2007), 269--279.


\bibitem{Burdick} B.S. Burdick, {\em A note on
completeness of hyperspaces}, in: General Topology
and Applications: Fifth Northeast Conference (Susan
Andima et al., ed.), Dekker (1991), pp. 19--24.

\bibitem{Burdick2} B.S. Burdick, {\em Representing quasi-uniform
spaces as the primes of ordered spaces}, preprint (2008).


\bibitem{CaoReillyRomaguera} J. Cao, I.L. Reilly
and S. Romaguera, {\em Some properties of quasi-uniform
multifunction spaces}, J. Austral. Math. Soc. (Ser.
A) {\bf 64} (1998), 169--177.


\bibitem{Deak} J. De\'ak, {\em A bitopological view of quasi-uniform
completeness. II}, Studia Sci. Math. Hung. {\bf 30} (1995), 411--431.


\bibitem{Doitchinov} D. Doitchinov, {\em On completeness of quasi-uniform
spaces}, C.R. Acad. Bulg. Sci. (7){\bf 41} (1988), 5--8.


\bibitem{Doitchinovstable} D. Doitchinov, {\em $E$-completions of quasi-uniform spaces},
Symposium on Categorical Topology (Rondebosch, 1994), 89--102, Univ. Cape Town,
Rondebosch, 1999.


\bibitem{Engelking} R. Engelking,
{\em General Topology}, Heldermann, 1989.

\bibitem{FedorchukKunzi}
V.V. Fedorchuk and H.-P.A. K\"unzi,
{\em Uniformly open mappings and uniform embeddings of function spaces},
Topology Appl. {\bf 61} (1995), 61--84.


\bibitem{FletcherLindgren} P. Fletcher and W.F. Lindgren, {\em
Quasi-uniform Spaces}, Dekker, New York, 1982.

\bibitem{Hohti} A. Hohti, {\em On supercomplete uniform spaces. V. Tamano's
product problem}, Fund. Math. {\bf 136} (1990), 121--125.

\bibitem{Isbell} J.R. Isbell, {\em Supercomplete spaces},
Pacific J. Math. {\bf 12}(1) (1962), 287--290.

\bibitem{Kunzi} H.-P.A. K\"unzi, {\em Nonsymmetric topology},
Topology with applications (Szeksz\'ard, 1993), Bolyai Soc. Math. Stud. 4.
Bolyai Math. Soc., Budapest 1995, pp. 303--338.


\bibitem{KunziIntroduction} H.-P.A. K\"unzi, {\em An introduction
to quasi-uniform spaces}, Chapter in: Beyond Topology, eds. F.
Mynard and E. Pearl, Contemporary Mathematics, Amer. Math. Soc.,
to appear.

\bibitem{KunziJunnila} H.-P. A. K\"unzi and H.J.K. Junnila, {\em
Stability in quasi-uniform spaces and the inverse problem},
Topology Appl. {\bf 49} (1993), 175--189.

\bibitem{KunziMrsevicReillyVamanamurthy} H.-P. A. K\"unzi,
M. Mr\v sevi\'c, I.L. Reilly and M.K. Vamanamurthy, {\em
Convergence, precompactness and symmetry in quasi-uniform
spaces}, Math. Japonica {\bf 38} (1993), 239--253.


\bibitem{KunziRomaguera} H.-P. A. K\"unzi and S. Romaguera,
{\em Well-quasi-ordering and the Hausdorff quasi-uniformity},
Topology Appl. {\bf 85} (1998), 207--218.


\bibitem{KunziRyser}
H.-P.A. K\"unzi and C. Ryser,
{\em The Bourbaki quasi-uniformity}, Topology Proc.
 {\bf 20}
(1995), 161--183.





\bibitem{RodriguezRomaguera} J. Rodr\'{\i}guez-L\'opez and
S. Romaguera, {\em The relationship between the Vietoris
topology and the Hausdorff quasi-uniformity}, Topology
Appl. {\bf 124} (2002), 451--464.

\bibitem{RomagueraSanchis} S. Romaguera and M. Sanchis, {\em
Completeness of hyperspaces on topological groups},
Journ. Pure Appl. Algebra {\bf 149} (2000), 287--293.

\bibitem{SanchezGranero} M.A. S\'anchez-Granero, {\em Weak
completeness of the Bourbaki quasi-uniformity},
Appl. Gen. Topology {\bf 2} (2001), 101--112.

\end{thebibliography}
\end{document}